\documentclass[12pt]{article}
\usepackage[utf8]{inputenc}

\usepackage{cite} 
\usepackage{amssymb} 
\usepackage{amsmath}
\usepackage{graphicx}
\usepackage{hyperref}
\usepackage{fullpage}

\long\def\beginpgfgraphicnamed#1#2\endpgfgraphicnamed{\includegraphics{#1}}

\usepackage{subfigure}

\renewcommand{\fbox}{} 

\newcommand{\R}{\mathbb{R}}
\newcommand{\Order}{\mathcal{O}}

\newtheorem{definition}{Definition}
\newtheorem{theorem}{Theorem}
\newcommand{\proof}{{\noindent\bf{Proof.}}}
\newcommand{\QED}{\hfill$\square$\bigbreak}

\begin{document} 

\title{Robust heteroclinic cycles in pluridimensions}
\author{Sofia B. S. D. Castro\footnote{Centro de Matem\'atica and Faculdade de Economia, Universidade do Porto, Portugal; email:\href{mailto:sdcastro@fep.up.pt}{sdcastro@fep.up.pt}}\ \ and Alastair M. Rucklidge\footnote{School of Mathematics, University of Leeds, Leeds LS2 9JT, UK; email: \href{mailto:A.M.Rucklidge@leeds.ac.uk}{A.M.Rucklidge@leeds.ac.uk}}}
\date{}

\maketitle

\begin{abstract}
\noindent
Heteroclinic cycles are sequences of equilibria along with trajectories that connect them in a cyclic manner.
We investigate a class of robust heteroclinic cycles that does not satisfy the usual condition that all connections between equilibria lie in flow-invariant subspaces of equal dimension.
We refer to these as robust heteroclinic cycles in pluridimensions.
The stability of these cycles cannot be expressed in terms of ratios of contracting and expanding eigenvalues in the usual way because, when the subspace dimensions increase, the equilibria fail to have contracting eigenvalues.
We develop the stability theory for robust heteroclinic cycles in pluridimensions, allowing for the absence of contracting eigenvalues.
We present four new examples, each with four equilibria and living in four dimensions, that illustrate the stability calculations.
Potential applications include modelling the dynamics of evolving populations when there are transitions between equilibria corresponding to mixed populations with different numbers of species.

\vspace{1ex}


\vspace{1ex}


\end{abstract}

\section{Introduction}\label{sec:introduction}

There has been a long-standing interest in robust heteroclinic cycles (sequences of equilibria and connecting trajectories between two consecutive equilibria in a cyclic manner) dating back to dos Reis~\cite{dosReis1984} and Guckenheimer and Holmes~\cite{Guckenheimer1988}. 
Robustness is achieved by ensuring that the connection between two saddle equilibria is of saddle-sink type in a lower-dimensional flow-invariant subspace. 
Flow-invariance appears naturally in problems with symmetry and in applications to population dynamics and game theory.
In the former, the flow-invariant subspaces are fixed-point subspaces, that is, sets of points that are preserved under the action of a symmetry subgroup. 
In the latter, extinction hyperplanes, where one or more variables is equal to zero, and subspaces of points with non-negative coordinates, provide the flow-invariance that leads to robust heteroclinic cycles.

A natural question to ask is whether or not robust heteroclinic cycles are stable, that is, are trajectories that start close to a cycle attracted closer to that cycle.
The systematic investigation of this started with Krupa and Melbourne \cite{Krupa1995,Krupa2004}, who established necessary and sufficient conditions for stability of certain types of ``simple'' cycles.
Their definition of simple (given below) includes the requirement that the flow-invariant subspaces containing the heteroclinic connections are all of the same dimension (and all of dimension two).
The condition for stability is expressed in terms of products and ratios of three types of eigenvalues of the Jacobian matrix at each of the equilibria.
The eigenvalues are called contracting, expanding and transverse (defined below), and the theory relies on the presence of both contracting and expanding eigenvalues, with transverse eigenvalues being optional.

There are a small number of published examples of heteroclinic cycles that do not fit this standard view~\cite{Sikder1994, Hawker2005a, Matthews1996, Rucklidge1995b}.
In these examples, there are equilibria in the cycle where there are no contracting eigenvalues, and so the existing theory for computing stability fails, although stability can be established in an \textit{ad hoc} fashion.
The reason for the absence of contracting eigenvalues is that the dimensions of the flow-invariant subspaces that contain the heteroclinic connections vary around the cycle (see Theorem~\ref{th:pluri-implies-no-contracting} below).

In this article we describe heteroclinic cycles that have at least two of the heteroclinic connections contained in flow-invariant subspaces of different dimension.
We call these \emph{robust heteroclinic cycles in pluridimensions}.
We focus on heteroclinic cycles (rather than networks), so we work with equilibria that all have one-dimensional unstable manifolds.
We take the first steps towards developing a more general stability theory, building on the work of Krupa and Melbourne~\cite{Krupa2004}.
We start by constructing four examples of robust heteroclinic cycles in pluridimensions.
Each of these examples has four equilibria and lives in~$\R^4$, and we show, with reasonable simplifying assumptions, that these are the four simplest examples.
We develop new stability results for these four examples, but we present the theory in a way that it can be readily applied to heteroclinic cycles in higher dimensions.
Furthermore, we envisage applications of our new approach to (for instance) multi-species Lotka--Volterra systems~\cite[Chapter~5]{Hofbauer1998} and heteroclinic networks, in particular, depth-two heteroclinic networks~\cite{Chawanya1997}.

Section~\ref{sec:preliminaries} recalls existing concepts relevant to our study. 
Section~\ref{sec:pluri} defines the object of our study: robust heteroclinic cycles in pluridimensions.
We show that having heteroclinic connections in flow-invariant subspaces with different dimensions around the cycle implies the absence of contracting eigenvalues at at least one equilibrium point.
With some reasonable assumptions, we present the four simplest examples of robust heteroclinic cycles in pluridimensions.
In Section~\ref{sec:return_maps}, we provide all the information required to calculate the return maps around these four examples and so determine their stability.
Although we present our results for these specific examples in~$\R^4$, we do so in a general way that can be extended to other examples in higher dimensions.
Section~\ref{sec:examples} presents numerical examples of each of the four robust heteroclinic cycles in pluridimensions given in Section~\ref{sec:pluri}, illustrating the dynamics for parameter values on both sides of the stability boundary.
We summarise and present ideas for future work in Section~\ref{sec:discussion}.

\section{Preliminaries}\label{sec:preliminaries}

Given an ordinary differential equation (ODE)
$$
\dot{x} = f(x),
$$
where $x \in \R^n$ and $f \in C^{\infty}$, we say that there is a heteroclinic cycle if there exist finitely many equilibria and trajectories connecting them in a  unique cyclic way.
Although saddle-saddle heteroclinic connections are in general not robust, when the ODE has flow-invariant subspaces, connections within these subspaces can occur as saddle-sink connections and are then robust.

Robust heteroclinic cycles have been classified as \emph{simple} in at least two distinct ways.
The first is given in the context of symmetric dynamics, where the ODE is equivariant under a group~$\Gamma$, and the flow-invariant spaces containing the connections are fixed-point subspaces of subgroups of~$\Gamma$.
The standard notation~\cite{Krupa2004} is to denote the invariant subspace containing the connection from an equilibrium $\xi_j$ to an equilibrium $\xi_{j+1}$, $[\xi_j \rightarrow \xi_{j+1}]$, by $P_j=\mbox{Fix}(\Sigma_j)$, where $\Sigma_j$ is an isotropy subgroup of~$\Gamma$, and to define $L_j=P_{j-1}\cap P_j$, so $\xi_j\in L_j$.
\begin{definition}[Krupa and Melbourne \cite{Krupa2004}]\label{def:simple-KM}
A robust heteroclinic cycle $X$ is \emph{simple} if $\dim P_j=2$ for all $j$ and $X$ intersects at most once each connected component of $L_j\backslash \{0\}$. 
\end{definition}
The second definition of simple applies even when the flow-invariant subspaces do not arise because of symmetry.
\begin{definition}[Hofbauer and Sigmund \cite{Hofbauer1998}]\label{def:simple-HS}
A robust heteroclinic cycle $X$ is \emph{simple} if the unstable manifold of every equilibrium $\xi_j \in X$, $W^u(\xi_j)$, has dimension 1.
\end{definition}
These two definitions are disjoint.
A heteroclinic cycle can satisfy Definition~\ref{def:simple-KM} but not Definition~\ref{def:simple-HS} if the connections are contained in 2-dimensional fixed-point spaces but there is more than one positive eigenvalue at an equilibrium.
On the other hand, examples of cycles that satisfy Definition~\ref{def:simple-HS} but not Definition~\ref{def:simple-KM} can be found in~$\R^3$ in Hawker and Ashwin~\cite{Hawker2005a} and in~$\R^4$ in Sikder and Roy~\cite{Sikder1994}.
In the example from~\cite{Hawker2005a}, the heteroclinic connections are contained in coordinate axes ($\dim P_j=1$) and in coordinate planes, and
in the example from~\cite{Sikder1994}, the connections are contained in planes and in three-dimensional spaces ($\dim P_j=3$).

Following the notation of \cite{Krupa2004}, we represent by~$P_j$ the smallest possible flow-invariant subspace containing the trajectory connecting $\xi_j$ to $\xi_{j+1}$ and use $L_j=P_{j-1}\cap P_j$ as above.
Note that this differs from the definition in~\cite{Krupa2004}, where the~$P_j$ are defined as subspaces fixed by a symmetry subgroup.

We classify the eigenvalues of the Jacobian matrix at an equilibrium $\xi_j$ as in \cite{Krupa2004}, that is, an eigenvalue is
\begin{itemize}
    \item \emph{radial} ($r$) if the corresponding eigenvector belongs to $L_j$;
    \item \emph{contracting} ($c$) if the corresponding eigenvector belongs to $P_{j-1}$ but not to $L_j$; call the space spanned by these eigenvectors $V_j(c)$;
    
    \item \emph{expanding} ($e$) if the corresponding eigenvector belongs to $P_{j}$ but not to $L_j$; call the space spanned by these eigenvectors $V_j(e)$;
    
    \item \emph{transverse} ($t$), all remaining eigenvalues. 
\end{itemize}
As explained in the Introduction, we focus here on heteroclinic cycles, not networks, and in particular we focus on heteroclinic cycles that contain only equilibria with one-dimensional unstable manifolds.
In particular, each equilibrium has a positive expanding eigenvalue and all other eigenvalues (radial, contracting and transverse) are negative.
We also simplify the presentation by considering at this stage only systems that have a single (positive) expanding eigenvalue, and avoid the complication of having a mixture of positive and negative expanding eigenvalues, as in the examples of Matthews and co-authors~\cite{Rucklidge1995b,Matthews1996}.

As well as Definitions~\ref{def:simple-KM} and~\ref{def:simple-HS}, there are several further classifications of heteroclinic cycles: simple of types~A, B and~C~\cite{Krupa2004}, simple of type~Z~\cite{Podvigina2012}, pseudo-simple~\cite{Podvigina2017} and quasi-simple~\cite{GarridodaSilva2019}. We do not need the details of these definitions here, as what we are about to introduce is different from all of these.
In particular, we depart from the assumption that $\dim P_j=2$ for all $j$ in Definition~\ref{def:simple-KM}, and from the even weaker assumption that $\dim P_j=\dim P_{j+1}$ for all $j$, which is used for type~$Z$ cycles by Podvigina~\cite{Podvigina2012} and for quasi-simple cycles by Garrido-da-Silva and Castro~\cite{GarridodaSilva2019}.
The fact that we do not rely on the equivariance of the vector field avoids some issues in determining the best definition of simple identified by the authors of~\cite{Podvigina2015a}, \cite{Podvigina2017}, and \cite{Chossat2018}.

\section{Robust cycles in pluridimensions}\label{sec:pluri}

We focus on problems where the flow-invariant connecting subspaces~$P_j$ do not all have the same dimensions, and define:
\begin{definition}\label{def-pluridimensional}
A robust heteroclinic cycle~$X$ is said to be a \emph{robust cycle in pluridimensions} if there exist two flow-invariant connecting subspaces with different dimensions, that is, $\dim P_{j-1}\neq \dim P_{j}$ for some~$j$.
\end{definition}
The examples from~\cite{Sikder1994,Hawker2005a,Castro2023,Matthews1996,Rucklidge1995b} all fit this definition, and many more examples can be constructed.
In order to give context to the stability theory developed in Section~\ref{sec:return_maps}, we aim in this section to generate the simplest possible examples of robust heteroclinic cycles in pluridimensions, and so we make a number of simplifying assumptions that nonetheless capture the essential features of cycles in pluridimensions, working in spaces of lowest dimension possible.
In Section~\ref{sec:discussion}, we explain how the stability theory applies, or can be extended, in the case of cycles in pluridimensions that do not satisfy each of these assumptions.

The simplifying assumptions that we make for generating examples are:
\begin{enumerate}
    \item[(A1)] $\dim W^u(\xi_i)=1$ for all equilibria $\xi_i \in X$;
    \item[(A2)] all coordinate axes and hyperplanes are flow-invariant subspaces;
    \item[(A3)] there is at most one equilibrium per connected component\footnote{By \emph{connected component}, we follow the meaning in~\cite{Krupa2004}: a coordinate axis has two connected components separated by the origin, a coordinate plane has four, separated by the coordinate axes, etc.} of each flow-invariant subspace;
    \item[(A4)] the origin is not part of the heteroclinic cycle.
\end{enumerate}
With Assumption~(A1), that all unstable manifolds of equilibria in the cycle are one dimensional, we remain within simple heteroclinic cycles by Definition~\ref{def:simple-HS}.
Assumption~(A2) is natural in the context of population dynamics and game theory, and avoids some of the complexities that arise in the presence of symmetry.
This assumption implies that the variables cannot cross coordinate planes and so cannot change sign, and that all eigenspaces are flow invariant.
Assumptions~(A1) and~(A2) together make the order in which trajectories visit the equilibria respect the order chosen for numbering the coordinates.
These assumptions also enable the stability calculations in Section~\ref{sec:return_maps}.
Assumptions (A3) and~(A4) come from~\cite{Krupa2004} and reduce the number of possibilities we have to consider in developing the classification in this section.
A weaker version of~(A3), as made by~\cite{Krupa2004}, would permit more that one equilibrium on (for example) a positive coordinate axis, but only one of these would be part of the cycle.
We make the stronger assumption~(A3) in order to reduce the number of possible equilibria and the range of possible heteroclinic cycles in pluridimensions.

The examples of heteroclinic cycles in pluridimensions in the literature~\cite{Castro2023, Matthews1996,Rucklidge1995b, Sikder1994, Hawker2005a} each have features that complicate the presentation of a general theory, and we use these simplifying assumptions to develop examples that illustrate the theory without additional complications.
The example in~\cite{Castro2023} has one equilibrium with a two-dimensional unstable manifold and so does not satisfy~(A1), and allows trajectories to leave one of the equilibria in a range of different directions.
In the convection and magnetoconvection examples in~\cite{Matthews1996,Rucklidge1995b} not all coordinate axes are flow invariant and so these do not satisfy~(A2). In these examples, there are equilibria with negative expanding eigenvalues, so the expanding directions are higher dimensional than strictly necessary.
The example in~\cite{Sikder1994} has two equilibria in a coordinate plane and so does not satisfy~(A3).
This is not a significant issue from the point of view of the stability theory, but in terms of generating new examples, allowing multiple equilibria in a coordinate axis or coordinate plane would lead to examples with arbitrarily many equilibria.
Finally, the example in~\cite{Hawker2005a} includes the origin in the heteroclinic cycle and so does not satisfy~(A4).
Again, this is not a significant issue from the point of view of the stability theory.

For the remainder of this section, we use the index~$j$ to refer to equilibria that satisfy $\dim P_{j-1}\neq \dim P_{j}$, and the index~$i$ for any equilibrium point, with no restriction on $\dim P_{i-1}$ and~$\dim P_{i}$.

In seeking examples of heteroclinic cycles in pluridimensions, there are several consequences of these assumptions, stated below.
\begin{itemize}
\item[(C1)]
It follows from (A1) that all radial, contracting and transverse eigenvalues are negative, otherwise there would be an equilibrium with an unstable manifold of dimension higher than one.

\item[(C2)]
It follows from (A1), (A2) and (A3) that $\dim P_i=\dim L_i+1$ for all $i$, and in addition from (A4) that $\dim L_i\geq1$.
This is because $P_{i}$ is the space containing the connection $[\xi_i \rightarrow \xi_{i+1}]$, and this connection is one-dimensional from~(A1) and is not in~$L_i$ from~(A3).
There can only be a single (positive) expanding eigenvalue from~(A2), so leaving $\xi_i$ can only increase the dimension by~$1$.
Assumption~(A4) gives us that $\dim L_i\geq1$.

\item[(C3)]
Since our definition of the flow-invariant connecting subspaces~$P_i$ is that they are the smallest subspaces that contain each connection $[\xi_i \rightarrow \xi_{i+1}]$, the number of non-zero coordinates on each connection is equal to the dimension of~$P_i$.
Similarly, since $\xi_i\in L_i=P_i\cap P_{i-1}$, the dimension of $L_i$ is equal to the number of non-zero coordinates of~$\xi_i$. 

\item[(C4)]
There is an~$i$ such that $\dim L_{i-1}\neq \dim L_{i}$.
This follows from~(C2) ($\dim P_i=\dim L_i+1$) and from the definition of cycles in pluridimensions ($P_{j}$ and $P_{j-1}$ have different dimensions).
Then there is at least one equilibrium point that is not on a coordinate axis.

\end{itemize}
Note that (C4) does not hold unless the cycle is in pluridimensions.

Standard results on the stability of robust heteroclinic cycles~\cite{Krupa1995} rely on all equilibria in the cycle having contracting and expanding directions.
However, a feature of robust cycles in pluridimensions is that some equilibria do not have contracting directions, and so standard stability results cannot be applied.
In addition, some equilibria have more than one contracting direction.
We show this in the following.
\begin{theorem}\label{th:pluri-implies-no-contracting}
For a robust cycle in pluridimensions~$X$ satisfying (A1), there is at least one equilibrium whose Jacobian matrix does not have contracting eigenvalues, and there is at least one equilibrium whose Jacobian matrix has at least two contracting eigenvalues.
\end{theorem}

\proof{}
For the first part of the theorem, note that if $L_i=P_{i-1}$ then there are no contracting eigenvalues at~$\xi_i$, since the contracting direction $V_i(c)$ is the empty set.
Because $X$ is a robust cycle in pluridimensions, we know that for some $j$ we have $\dim P_{j} > \dim P_{j-1}$.
Then $\dim P_j \ominus L_j > \dim P_{j-1} \ominus L_j$, where we use $A\ominus B$ to denote the orthogonal complement of set~$B$ inside set~$A$.
The set $P_j \ominus L_j$ is the expanding direction~$V_j(e)$, which is one dimensional by Assumption~(A1).
So $1>\dim P_{j-1} \ominus L_j$, so $\dim P_{j-1} \ominus L_j=0$, and the contracting direction~$V_j(c)$ is empty and so for this equilibrium point, there are no contracting eigenvalues.

For the second part of the theorem, we know that for some $j$ we have $\dim P_{j-1} > \dim P_{j}$.
Then $\dim P_{j-1} \ominus L_{j} > \dim P_{j} \ominus L_{j}$.
The set $P_{j} \ominus L_{j}$ is the expanding direction~$V_j(e)$, which is one dimensional by Assumption~(A1).
So $\dim P_{j-1} \ominus L_{j}\geq 2$, and the contracting direction~$V_j(c)$ is at least two-dimensional. Hence, for this equilibrium point, there are at least two contracting eigenvalues.
\QED

We note that Theorem~\ref{th:pluri-implies-no-contracting} does not hold for robust cycles that are not pluridimensional: such cycles have contracting eigenvalues at all equilibria.
Conversely, the example of~\cite{Castro2022b} is pluridimensional according to Definition~\ref{def-pluridimensional}, but does not satisfy~(A1), and all equilibria in that cycle do have contracting dimensions.

The second part of Theorem~\ref{th:pluri-implies-no-contracting} highlights another feature of robust cycles in pluridimensions, which has consequences for the calculation of the stability of the cycle.
These are addressed in Subsection~\ref{sec:localmaps} and illustrated in Figure~\ref{fig:logh}.

We now construct examples of robust heteroclinic cycles in pluridimensions, aiming for examples that demonstrate the typical features of these cycles.
Our assumptions avoid the distraction of higher dimensional unstable manifolds, complications arising from symmetry considerations, having multiple routes in and out of invariant subspaces, and the special case of connections along coordinate axes.
We show in Theorems~\ref{th:no-cycle-2-equil}--\ref{th:no-cycle-3-equil} that robust
cycles in pluridimensions satisfying (A1)--(A4) have to be in at least four
dimensions and have to have at least four equilibria, so our examples will be 
in~$\R^4$ and will have four equilibria. 

\begin{theorem}\label{th:no-cycle-2-equil}
There are no robust cycles in pluridimensions satisfying (A1)--(A4) with exactly two equilibria.
\end{theorem}

\proof{}
We have $\xi_1\in L_1=P_1\cap P_2$, and $\xi_2\in L_2=P_2\cap P_1$, so $L_1=L_2$, and $\xi_1$ and $\xi_2$ are in the same space.
This does not contradict Assumption~(A3), since the two equilibria could be on a coordinate axis with the origin in between.
From consequence~(C2), which relies on (A1)--(A4), we have $\dim P_i=\dim L_i+1$, so $\dim P_1=\dim P_2$, and the cycle is not in pluridimensions.
\QED

\begin{theorem}\label{th:no-cycle-R3}
There are no robust cycles in pluridimensions satisfying (A1)--(A4) in $\R^3$.
\end{theorem}

\proof{}
The origin is not part of the cycle by~(A4), so the dimensions of the $L$ subspaces are at least one.
By~(C2), the dimensions of the $P$ subspaces are at least two.
The $P$ subspaces cannot be three-dimensional in~$\R^3$, since the equilibrium at the end of the connection in this three-dimensional space is a sink, and therefore it cannot have an unstable manifold in~$\R^3$.
Therefore all the $P$ subspaces are two dimensional, and the cycle is not in pluridimensions.
\QED

\begin{theorem}\label{th:no-cycle-3-equil}
There are no robust cycles in pluridimensions satisfying (A1)--(A4) in $\R^4$ with exactly three equilibria.
\end{theorem}

\proof{}
In $\R^4$, by an argument similar to the proof of Theorem~\ref{th:no-cycle-R3}, there are not enough dimensions for a pluridimensional cycle if the cycle does not include an equilibrium point on an axis. 
So, without loss of generality, let $\xi_1$ be on the $x_1$ axis and let $\xi_2$ be in the $(x_1,x_2)$ plane.
This means that $L_1=\{(x_1,0,0,0):x_1\in\R\}$ and $P_1=L_2=\{(x_1,x_2,0,0):x_1,x_2\in\R\}$. 
Then by~(C2), we must have $P_2=\{(x_1,x_2,x_3,0):x_1,x_2,x_3\in\R\}$.

The third equilibrium point $\xi_3\in P_2$, so it must have three, two or one non-zero coordinates.
The connection $\xi_3\rightarrow\xi_1$ must have $x_2=0$, because otherwise $\xi_1$ would not be a sink in the relevant subspace.
This implies that $\xi_3$ itself must have $x_2=0$.
In addition, $\xi_3$ must be unstable in the $x_4$ direction (else, Theorem~\ref{th:no-cycle-R3} would apply), and so in order to get back to~$\xi_1$, given~(A1), $\xi_3$~must have $x_1\neq0$.
By~(A3), $\xi_3\notin L_1$, so $\xi_3$ must have $x_1\neq0$ and $x_3\neq0$.
This means that $L_3=\{(x_1,0,x_3,0):x_1,x_3\in\R\}$, and $P_3=\{(x_1,0,x_3,x_4):x_1,x_3,x_4\in\R\}$.

The implication of this is that $\xi_1$ and $\xi_3$ are both sinks in~$L_3$, a two-dimensional space with $\xi_1$ on the $x_1$~axis and $\xi_3$ in the $(x_1,x_3)$ plane.
There are no additional equilibria in $L_3$ by~(A3), so there must be a periodic orbit surrounding~$\xi_3$ in the $(x_1,x_3)$ plane.
However, this periodic orbit has an invariant manifold associated with perturbations in the~$x_2$ direction.
This cylindrical manifold surrounds the heteroclinic connection from $\xi_2$ to $\xi_3$. 
The connection from $\xi_1$ to $\xi_2$ must therefore cross this manifold, which is a contradiction.
\QED

\begin{table}
\begin{center}
\begin{tabular}{|c|c|c|c|c|l|}
\hline
      & $\dim P_1$ & $\dim P_2$ & $\dim P_3$ & $\dim P_4$ & Illustrated in\\
     \hline 
    Case 1 & 2 & 3 & 2 & 3 & Figure~\ref{fig:case1}\\
    \hline
    Case 2 & 2 & 3 & 2 & 2 & Figure~\ref{fig:case2}\\
    \hline
    Case 3 &2 & 3 & 3 & 2 & Figure~\ref{fig:case3}\\
    \hline
    Case 4 &2 & 3 & 3 & 3 & Figure~\ref{fig:case4}\\
    \hline
\end{tabular}
\end{center}
\caption{Four cases of heteroclinic cycles in pluridimensions in~$\R^4$ satifying (A1)--(A4).
The dimensions of the $L_i$ subspaces satisfy $\dim L_i=\dim P_i-1$ by~(C2).
\label{table:fourcases}}
\end{table}

\begin{figure}
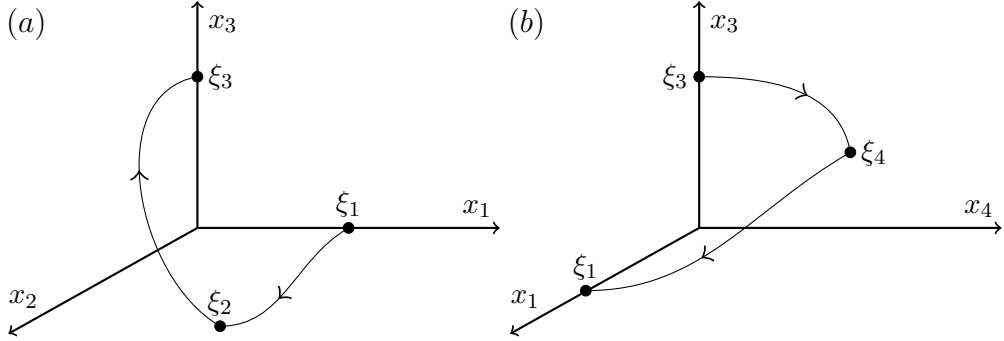

\centering
\subfigure{\fbox{%
\beginpgfgraphicnamed{Case_1_xi1_xi2_xi3}%
\endpgfgraphicnamed}}
\subfigure{\fbox{%
\beginpgfgraphicnamed{Case_1_xi3_xi4_xi1}%
\endpgfgraphicnamed}}
\caption{Robust heteroclinic cycles in pluridimensions in $\R^4$: Case~1. 
The equilibria $\xi_1$ and $\xi_3$ are on axes, while $\xi_2$ and $\xi_4$ are in the $(x_1,x_2)$ and $(x_3,x_4)$ planes respectively.
The equilibria and connections $\xi_1 \rightarrow \xi_2 \rightarrow \xi_3$ are identical to those in Case~2.
For ease of presentation, in this and the next three figures, we split the four-dimensional space into two or more panels.
In this and all other figures, the arrows on trajectories indicate the direction of travel and arrows on coordinate axes indicate their orientation.
\label{fig:case1}}
\end{figure}

\begin{figure}
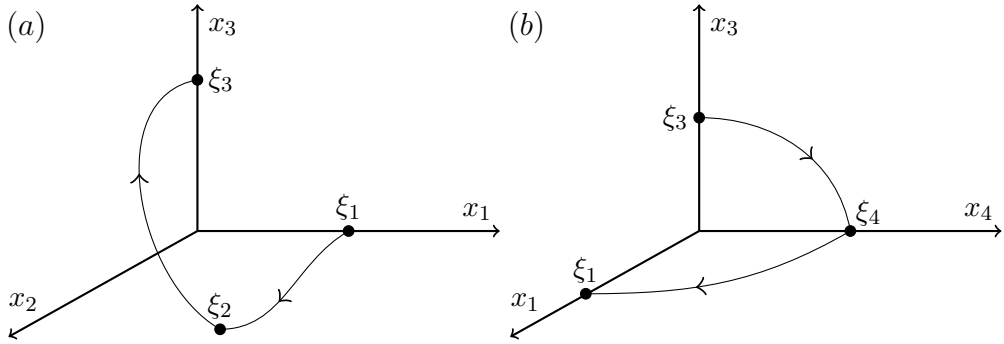

\centering
\subfigure{\fbox{%
\beginpgfgraphicnamed{Case_2_xi1_xi2_xi3}%
\endpgfgraphicnamed}}
\subfigure{\fbox{%
\beginpgfgraphicnamed{Case_2_xi3_xi4_xi1}%
\endpgfgraphicnamed}}
\caption{Robust heteroclinic cycles in pluridimensions in $\R^4$: Case~2.
The equilibria $\xi_1$, $\xi_3$ and $\xi_4$ are on axes, while $\xi_2$ is in the $(x_1,x_2)$ plane.
The equilibria and connections $\xi_1 \rightarrow \xi_2 \rightarrow \xi_3$ are identical to those in Case~1.
\label{fig:case2}}
\end{figure}

\begin{figure}
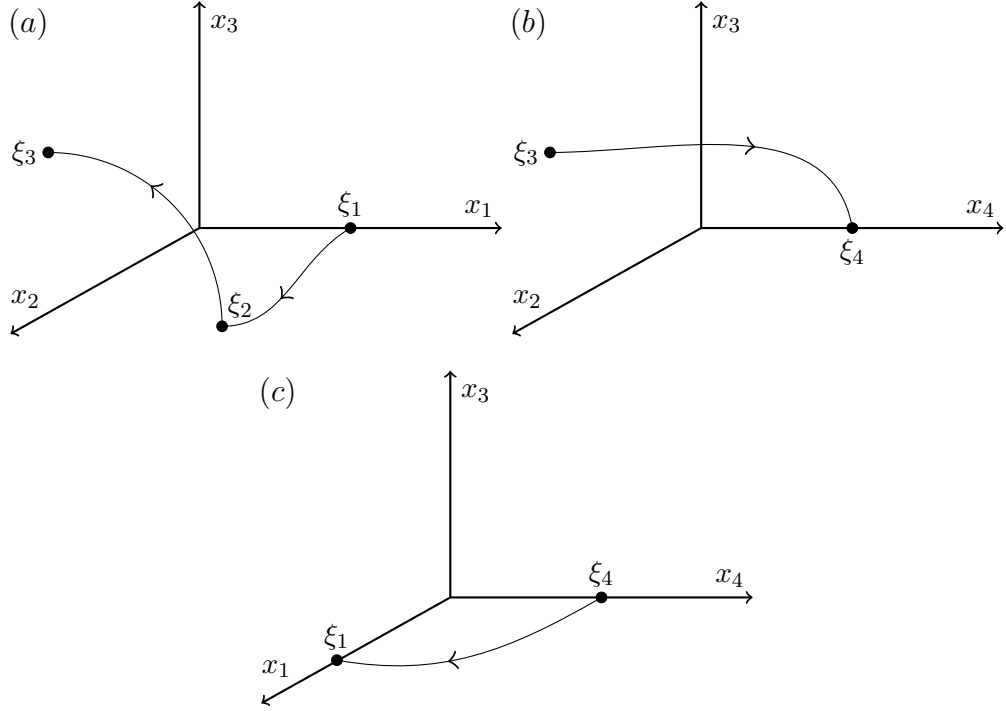

\centering
\subfigure{\fbox{%
\beginpgfgraphicnamed{Case_3_xi1_xi2_xi3}%
\endpgfgraphicnamed}}
\subfigure{\fbox{%
\beginpgfgraphicnamed{Case_3_xi3_xi4}%
\endpgfgraphicnamed}}
\subfigure{\fbox{%
\beginpgfgraphicnamed{Case_3_xi4_xi1}%
\endpgfgraphicnamed}}
\caption{Robust heteroclinic cycles in pluridimensions in $\R^4$: Case~3. 
The equilibria $\xi_1$ and $\xi_4$ are on axes, while $\xi_2$ and $\xi_3$ are in the $(x_1,x_2)$ and $(x_2,x_3)$ planes respectively.
The equilibria and connections $\xi_1 \rightarrow \xi_2 \rightarrow \xi_3$ are identical to those in Case~4.
\label{fig:case3}}
\end{figure}

\begin{figure}
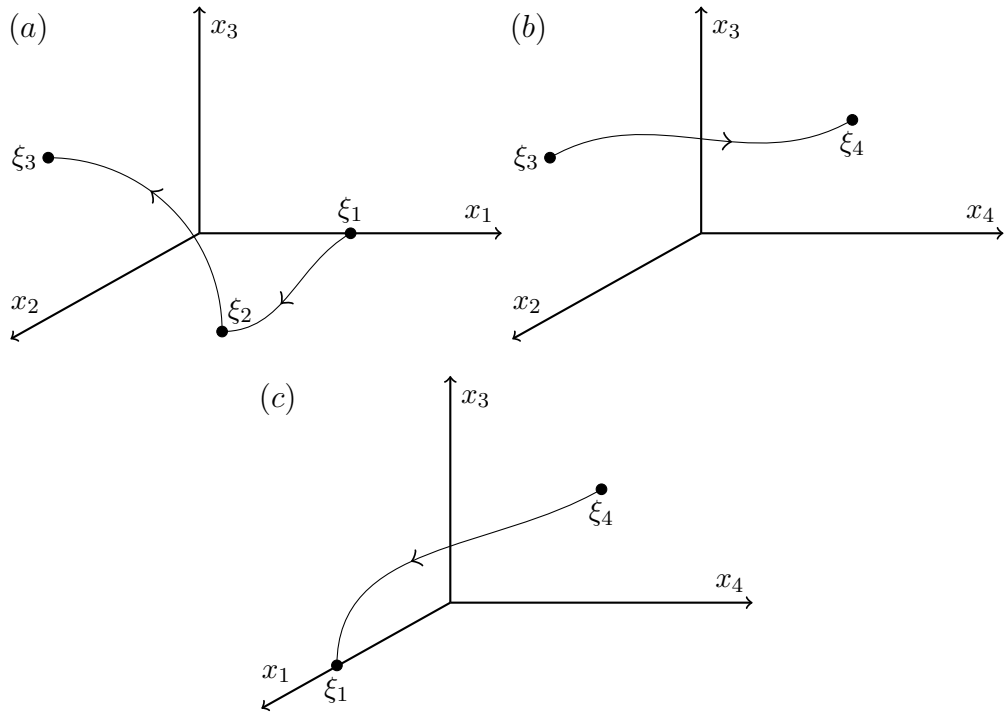

\centering
\subfigure{\fbox{%
\beginpgfgraphicnamed{Case_4_xi1_xi2_xi3}%
\endpgfgraphicnamed}}
\subfigure{\fbox{%
\beginpgfgraphicnamed{Case_4_xi3_xi4}%
\endpgfgraphicnamed}}
\subfigure{\fbox{%
\beginpgfgraphicnamed{Case_4_xi4_xi1}%
\endpgfgraphicnamed}}
\caption{Robust heteroclinic cycles in pluridimensions in $\R^4$: Case~4. 
The equilibria $\xi_1$ is on an axis, while $\xi_2$, $\xi_3$ and $\xi_4$ are in the $(x_1,x_2)$, $(x_2,x_3)$ and $(x_3,x_4)$ planes respectively.
The equilibria and connections $\xi_1 \rightarrow \xi_2 \rightarrow \xi_3$ are identical to those in Case~3.
\label{fig:case4}}
\end{figure}

We end this section by showing that there are only four heteroclinic cycles connecting four equilibria in pluridimensions in~$\R^4$ satisfying (A1)--(A4), apart from coordinate permutations. 
The $P$ subspaces cannot be one dimensional, from~(C2), and the $P$ subspaces cannot be four dimensional, since the equilibria at the end of a connection must be a sink within that subspace, and in $\R^4$ there would be no dimension for its unstable manifold.
In pluridimensions, the $P$ subspaces do not all have the same dimension, and there must be a mixture of $\dim P=2$ and $\dim P=3$.
Without loss of generality, we take $\dim P_1=2$ and $\dim P_2=3$, and so $\dim L_1=1$ and $\dim L_2=2$ by~(C2).
Then, there are only four possible choices of dimension of subspaces~$P_3$ and $P_4$, listed in Table~\ref{table:fourcases}.
The four examples are illustrated in Figures~\ref{fig:case1}--\ref{fig:case4}.

\paragraph{Cases 1 and 2.} 
Assumption~(A2) ensures that it suffices to look at positive coordinates, so let $x_i^{(j)} \in \R_+$, where $i$ refers to the index of the coordinate and $j$ to that of the equilibrium point. 
The first three equilibria are, given~(C3) and~(C4),
\begin{itemize}

\item  $\xi_1 = (x^{(1)}_1,0,0,0)$ and $L_1= \{ (x_1,0,0,0): \; x_1 \in \R \}$, since equilibria are not the origin by~(A4);

\item  $\xi_2 = (x^{(2)}_1,x^{(2)}_2,0,0)$ and $L_2=P_1= \{ (x_1,x_2,0,0): \; x_1, x_2 \in \R \}$, by our assumption on the dimensions of $P_1$ and $P_2$ and by~(C2);

\item  $\xi_3 = (0,0,x^{(3)}_3,0)$, $P_2= \{ (x_1,x_2,x_3,0): \; x_1, x_2, x_3 \in \R \}$ and $L_3= \{ (0,0,x_3,0): \; x_3 \in \R \}$, because $\xi_3$ is in a one-dimensional subspace from $\dim P_3=2$ in these two cases.
\end{itemize}
There are two different choices for $\xi_4$ with the same $P_3$ but different $L_4$ and $P_4$:

\noindent
Case~1: $\xi_4 = (0,0,x^{(4)}_3,x^{(4)}_4)$, $L_4=P_3= \{ (0,0,x_3,x_4): \; x_3, x_4 \in \R \}$ and $P_4=\{ (x_1,0,x_3,x_4): \; x_1,x_3,x_4 \in \R \}$;

\noindent
Case~2: $\xi_4 = (0,0,0,x^{(4)}_4)$, $P_3=\{ (0,0,x_3,x_4): \; x_3, x_4 \in \R \}$ and $L_4= \{ (0,0,0,x_4): \; x_4 \in \R \}$ and $P_4=\{ (x_1,0,0,x_4): \; x_1,x_4 \in \R \}$.

These two examples of robust heteroclinic cycles in pluridimensions in~$\R^4$ satisfying (A1)--(A4) are illustrated in Figures~\ref{fig:case1}--\ref{fig:case2}.

\paragraph{Cases 3 and 4.} 
Again, let $x_i^{(j)} \in \R_+$. 
The first two equilibria are the same as in Cases~1 and 2, as are the spaces $L_1$, $P_1$, $L_2$, and $P_2$. The other two equilibria are, given~(C3),
\begin{itemize}
    \item  $\xi_3 = (0,x^{(3)}_2,x^{(3)}_3,0)$, $P_2= \{ (x_1,x_2,x_3,0): \; x_1, x_2, x_3 \in \R \}$ and $L_3= \{ (0,x_2,x_3,0): \; x_2, x_3 \in \R \}$, because $\xi_3$ is in a two-dimensional subspace from $\dim P_3=3$ in these two cases.
\end{itemize}
Again, there are two different choices for $\xi_4$ with the same $P_3$ but different $L_4$ and $P_4$:

\noindent
Case~3:  $\xi_4 = (0,0,0,x^{(4)}_4)$, $P_3=\{ (0,x_2,x_3,x_4): \; x_2, x_3, x_4 \in \R \}$ and $L_4= \{ (0,0,0,x_4): \; x_4 \in \R \}$ and $P_4=\{ (x_1,0,0,x_4): \; x_1,x_4 \in \R \}$;

\noindent
Case~4: $\xi_4 = (0,0,x^{(4)}_3,x^{(4)}_4)$, $L_4=P_3= \{ (0,x_2,x_3,x_4): \; x_2, x_3, x_4 \in \R \}$ and $P_4=\{ (x_1,0,x_3,x_4): \; x_1,x_3,x_4 \in \R \}$.

These two examples of robust heteroclinic cycles in pluridimensions in~$\R^4$ satisfying (A1)--(A4) are illustrated in Figures~\ref{fig:case3}--\ref{fig:case4}.

\section{Stability of cycles in pluridimensions}\label{sec:return_maps}

In this section, we compute return maps for cycles in pluridimensions, starting with definitions of the incoming and outgoing cross sections at the equilibria and the global maps between them.
We then define local maps at each equilibrium point: it is natural to express these as linear maps, matrices multiplying the logarithms of the coordinates.
The dynamics from one equilibrium point to the next is described by composing a local map at the first equilibrium point with the global map that leads to the next.
For trajectories that are very close to the cycle, having very small values of some coordinates with very large negative values of their logarithms, the composed map is dominated by a \emph{transition matrix} multiplying the logarithms of the coordinates~\cite{Field1991}.
Multiplying these transition matrices around the cycle gives the overall behaviour of trajectories very close to the cycle, and in particular gives the stability of the cycle~\cite{Krupa1995}.

In cycles that are not in pluridimensions, the form of the map from one equilibrium point to the next does not depend on the previous equilibrium point in the cycle, as the dimensions of the $P$ subspaces are all the same.
A new feature of cycles in pluridimensions is that the form of the map from one equilibrium point to the next can depend on the location of the previous equilibrium point in the cycle: a map from a point on an axis to a point in a plane will be different depending on whether the previous point was on an axis or in a plane. 
This is explained in detail in Subsection~\ref{sec:transitionmatrices}.
A consequence of this is that the transition matrices are not necessarily square.

In this section, we use generic coordinates, represented either by $(z_i,z_{i+1},z_{i+2},z_{i+3})\in\R^4$, or by $(z_{i-1},z_i,z_{i+1},z_{i+2})\in\R^4$, or by $(z_{i-2},z_{i-1},z_i,z_{i+1})\in\R^4$, depending on which is more convenient (the choice is clearly stated where needed).
Using these $z$'s rather than~$x$'s distinguishes these general coordinates, labelled by $i-2$, \dots, $i+3$ etc., from the specific ones in the examples in Sections~\ref{sec:pluri} and~\ref{sec:examples}, labelled by $1$, $2$, $3$ and~$4$.

\subsection{Global maps and cross sections}\label{sec:globalmaps}

In this subsection, we establish the possible global maps for heteroclinic cycles in pluridimensions in~$\R^4$.
The four cases from Table~\ref{table:fourcases}, illustrated in Figures~\ref{fig:case1}--\ref{fig:case4}, contain mixtures of four types of transition, from axis-to-axis, axis-to-plane, plane-to-axis and plane-to-plane.
These four types of transition are illustrated in Figures~\ref{fig:axis-to-axis}--\ref{fig:plane-to-plane} below, and the four global maps between the equilibria are derived in the paragraphs below.

Near~$\xi_i$, the outgoing cross section in the direction of~$\xi_{i+1}$ is $H_i^{out,i+1}$ (defined more precisely in each case below), and the incoming cross section near~$\xi_{i+1}$ from the direction of~$\xi_i$ is $H_{i+1}^{in,i}$.
In this context, all coordinates are non-negative.
The global maps go from $H_i^{out,i+1}$ to $H_{i+1}^{in,i}$, and are denoted by~$\Psi_{i \to i+1}$.

\begin{figure}
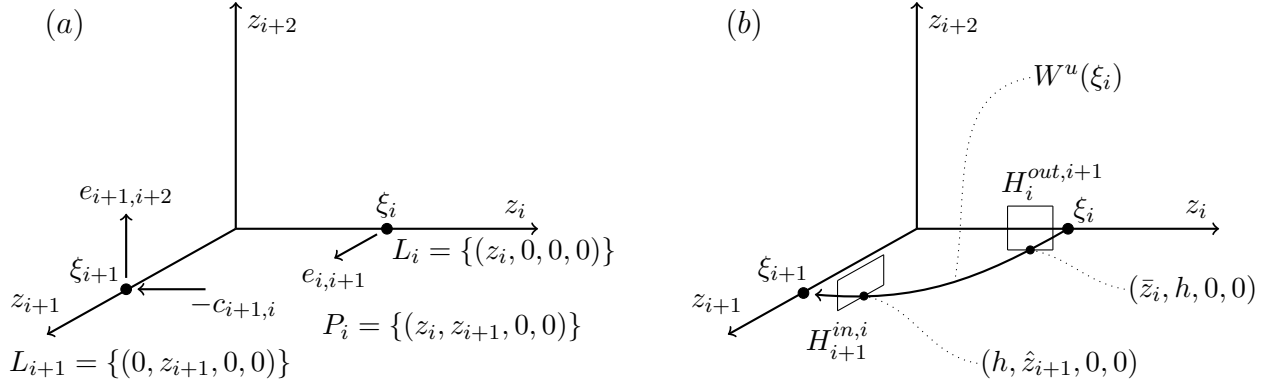

\hbox to \hsize{\fbox{%
\beginpgfgraphicnamed{Fig-global-map_AA_1}%
\endpgfgraphicnamed}%
\hfil\fbox{%
\beginpgfgraphicnamed{Fig-global-map_AA_2}%
\endpgfgraphicnamed}}
\caption{The global map between two equilibria on axes.
Panel~(a) shows the flow-invariant subspaces, equilibria, and the direction of the eigenvalues.
Panel~(b) shows the cross sections and coordinates of their intersection with the unstable manifold of $\xi_i$.
\label{fig:axis-to-axis}}
\end{figure}

\paragraph{Axis-to-axis connection.}
This type of global map describes the dynamics along a connection between two equilibria on coordinate axes, such as the connection between $\xi_3$ and $\xi_4$ in Case~2, Figure~\ref{fig:case2}(b).
Figure~\ref{fig:axis-to-axis} illustrates important points in the construction of the global map.

We define the cross sections as follows:
 $$
 H_i^{out,i+1}=\{ (z_i,z_{i+1},z_{i+2},z_{i+3}) \in \R^4: \;\; z_i=\Order(1), z_{i+1}=h, z_{i+2}<h, z_{i+3}<h \}
 $$
and
 $$
 H_{i+1}^{in,i}=\{ (z_i,z_{i+1},z_{i+2},z_{i+3}) \in \R^4: \;\; z_i=h, z_{i+1}=\Order(1), z_{i+2}<h, z_{i+3}<h \},
 $$
for small $h$. 
At~$\xi_i$, $z_i$ is the radial coordinate (and is of order~$1$) and $z_{i+1}$ is the expanding coordinate, and at~$\xi_{i+1}$, $z_{i+1}$~is the radial coordinate, $z_i$~is the contracting coordinate, and $z_{i+2}$~is the expanding coordinate.
The remaining coordinates can be contracting or transverse, depending on previous or subsequent connections.

The unstable manifold of $\xi_i$, $W^u(\xi_i)$, is the connection from $\xi_i$ to $\xi_{i+1}$ in $P_i$, where $z_{i+2}=z_{i+3}=0$.
It intersects the cross sections at
$$
W^u(\xi_i) \cap H_i^{out,i+1} = \{ (\bar{z}_i,h,0,0) \} 
\qquad\mbox{and}\qquad 
W^u(\xi_i) \cap H_{i+1}^{in,i} = \{ (h,\hat{z}_{i+1},0,0) \}.
$$
Throughout, we will use bar and hat accents to indicate the radial coordinate values where the unstable manifolds intersect outgoing and incoming sections respectively.
The plane $P_i$ is invariant so $W^u(\xi_i)$ leaves $\xi_i$ with $(z_{i+2},z_{i+3})=(0,0)$ and arrives at $\xi_{i+1}$ also with $(z_{i+2},z_{i+3})=(0,0)$.
We write a point in $H_i^{out,i+1}$ as $(z_i,z_{i+1},z_{i+2},z_{i+3})=(\bar{z}_i+\tilde{z}_i,h,z_{i+2},z_{i+3})$, with $|\tilde{z}_i|<h$.
Throughout, we will use tilde accents to indicate displacements from the unstable manifold in the radial direction at the outgoing section.
The linearisation around $W^u(\xi_i)$ provides the global map (using $[...]$ for entries in the matrix that are of no consequence for the study of stability)
\begin{equation} \label{eq:globalaxistoaxis}
    \Psi_{i \to i+1} \left(\begin{array}{c} 
z_i\\ z_{i+2}\\ z_{i+3}
    \end{array} \right)
    = \left(\begin{array}{ccc} 
[...] & [...] & [...] \\
0 & A^{i+2}_{i \to i+1} & 0 \\
0 & 0 & A^{i+3}_{i \to i+1}
    \end{array} \right) \left(\begin{array}{c} 
\tilde{z}_i\\ z_{i+2}\\ z_{i+3}
    \end{array} \right) + \left(\begin{array}{c} 
\hat{z}_{i+1}\\ 0\\ 0
    \end{array} \right).
\end{equation}
This map takes as argument the values of $(z_i,z_{i+2},z_{i+3})$ on the outgoing section and returns the values of $(z_{i+1},z_{i+2},z_{i+3})$ on the incoming section.
The invariance of the $z_{i+2}=0$ and the $z_{i+3}=0$ subspaces from~(A2) leads to the diagonal structure of the $(i+2,i+3)$ part of the matrix.
The $A$~coefficients come from the linearisation around the unstable manifold.
The un-named entries in the matrix, indicated by~$[...]$, lead to contributions that are small compared to the fixed order~$1$ value of~$\hat{z}_{i+1}$.
There are also $\Order(h^2)$ corrections to this linearised map (not written).
Both of these small contributions to the global map will be disregarded when all variables are small, close to the heteroclinic cycle.

\begin{figure}
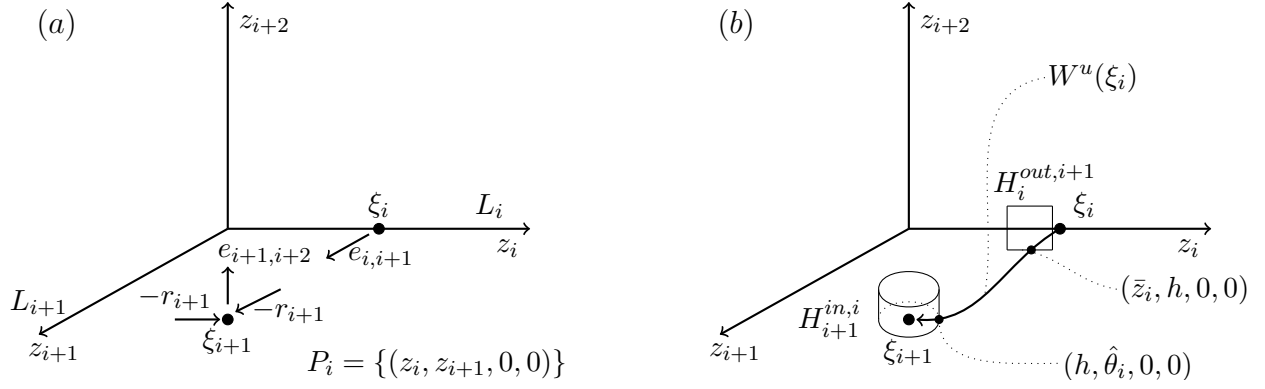

\hbox to \hsize{\fbox{%
\beginpgfgraphicnamed{Fig-global-map_AP_1}%
\endpgfgraphicnamed}%
\hfil\fbox{%
\beginpgfgraphicnamed{Fig-global-map_AP_2}%
\endpgfgraphicnamed}}
\caption{The global map between an equilibrium on an axis to an equilibrium in a plane. 
Panel~(a) shows the flow-invariant subspaces, equilibria, and the direction of the eigenvalues.
Panel~(b) shows the cross sections and coordinates of their intersection with the unstable manifold of $\xi_i$.
The two radial (possibly complex) eigenvalues at $\xi_{i+1}$ are both labelled~$-r_{i+1}$.
\label{fig:axis-to-plane}}
\end{figure}

\paragraph{Axis-to-plane connection.}
This type of global map describes the dynamics along a connection from an equilibrium on a coordinate axis to an equilibrium in a coordinate plane, containing that coordinate axis.
This situation is illustrated by the connection between $\xi_1$ and $\xi_2$ in Case~1, Figure~\ref{fig:case1}(a). 
In this case the connection is between an equilibrium with one non-zero coordinate and an equilibrium with two non-zero coordinates. 
Figure~\ref{fig:axis-to-plane} illustrates important points in the construction of the global map.

We define the outgoing cross section near $\xi_i$ as follows:
 $$
 H_i^{out,i+1}=\{ (z_i,z_{i+1},z_{i+2},z_{i+3}) \in \R^4: \;\; z_i=\Order(1), z_{i+1}=h, z_{i+2}<h, z_{i+3}<h \}
 $$
for small $h$. At~$\xi_i$, $z_i$ is the radial coordinate (and is of order~$1$) and $z_{i+1}$ is the expanding coordinate.
However, there is no contracting direction at~$\xi_{i+1}$, so the incoming section is defined in terms of the two radial coordinates, $(z_i,z_{i+1})$.
Because of this, it is natural to use polar coordinates in the radial direction centered on~$\xi_{i+1}$.
We define the displacement from $\xi_{i+1}$ as $(\rho_i\cos\theta_i,\rho_i\sin\theta_i)$, for $\rho_i\geq0$ and $0\leq\theta_i<2\pi$.
We define a cylinder for the incoming cross section near $\xi_{i+1}$ as follows:
 $$
 H_{i+1}^{in,i}=\{ (z_i,z_{i+1},z_{i+2},z_{i+3}) \in \R^4: \;\; \rho_i=h, 0\leq\theta_i<2\pi, z_{i+2}<h, z_{i+3}<h \},
 $$
for small $h$.

The unstable manifold of $\xi_i$, $W^u(\xi_i)$, which is the connection from $\xi_i$ to $\xi_{i+1}$ in~$P_i$, intersects the cross sections at
$$
W^u(\xi_i) \cap H_i^{out,i+1} = \{ (\bar{z}_i,h,0,0) \}
\qquad\mbox{and}\qquad 
W^u(\xi_i) \cap H_{i+1}^{in,i} = \{ (\rho_i=h,\theta_i=\hat{\theta}_{i},0,0) \}.
$$
The plane $P_i$ is invariant so $W^u(\xi_i)$ leaves $\xi_i$ with $(z_{i+2},z_{i+3})=(0,0)$ and arrives at $\xi_{i+1}$ also with $(z_{i+2},z_{i+3})=(0,0)$.
As before, we write a point in $H_i^{out,i+1}$ as $(z_i,z_{i+1},z_{i+2},z_{i+3})=(\bar{z}_i+\tilde{z}_i,h,z_{i+2},z_{i+3})$, with $|\tilde{z}_i|<h$.
The linearisation around $W^u(\xi_i)$ provides the global map
\begin{equation} \label{eq:globalaxistoplane}
    \Psi_{i \to i+1} \left(\begin{array}{c} 
{z}_{i}\\ z_{i+2}\\ z_{i+3}
    \end{array} \right)
    = \left(\begin{array}{ccc} 
[...] & [...] & [...] \\
0 & A^{i+2}_{i \to i+1} & 0 \\
0 & 0 & A^{i+3}_{i \to i+1}
    \end{array} \right) \left(\begin{array}{c} 
\tilde{z}_i\\ z_{i+2}\\ z_{i+3}
    \end{array} \right) + \left(\begin{array}{c} 
\hat{\theta}_{i}\\ 0\\ 0
    \end{array} \right).
\end{equation}
This map takes as argument the values of $(z_i,z_{i+2},z_{i+3})$ on the outgoing section and returns the values of $(\theta_i,z_{i+2},z_{i+3})$ on the incoming section.
The invariance of the $z_{i+2}=0$ and the $z_{i+3}=0$ subspaces leads to the diagonal structure of the $(i+2,i+3)$ part of the matrix.
The $A$~coefficients come from the linearisation around the unstable manifold.
The un-named entries in the matrix, indicated by~$[...]$, lead to contributions that are small compared to the fixed order~$1$ value of~$\hat{\theta}_i$.       
There are also $\Order(h^2)$ corrections to this linearised map (not written).
Both of these small contributions to the global map will be disregarded when all variables are small, close to the heteroclinic cycle.                                                                        

\begin{figure}
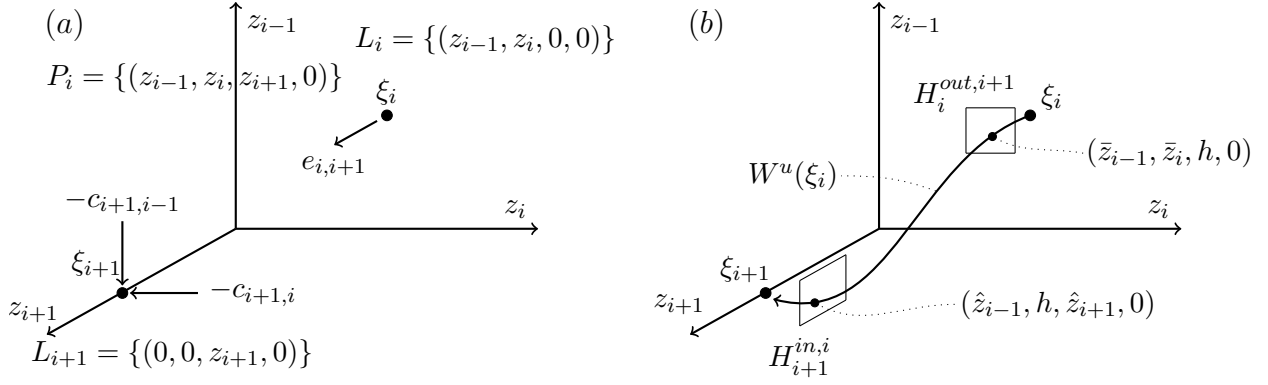

\hbox to \hsize{\fbox{%
\beginpgfgraphicnamed{Fig-global-map_PA_1}%
\endpgfgraphicnamed}%
\hfil\fbox{%
\beginpgfgraphicnamed{Fig-global-map_PA_2}%
\endpgfgraphicnamed}}
\caption{The global map between an equilibrium in a plane to an equilibrium on an axis.
Panel~(a) shows the flow-invariant subspaces, equilibria, and the direction of the eigenvalues.
Panel~(b) shows the cross sections and coordinates of their intersection with the unstable manifold of $\xi_i$.
We have illustrated the case where $z_i$ contracts more slowly than~$z_{i-1}$.
\label{fig:plane-to-axis}}
\end{figure}

\paragraph{Plane-to-axis connection.}
This type of global map describes the dynamics along a connection from an equilibrium in a coordinate plane to an equilibrium on a coordinate axis, not contained in that coordinate plane.
This situation is illustrated by the connection between $\xi_2$ and $\xi_3$ in Case~1, Figure~\ref{fig:case1}(a).
In this case, the connection is from an equilibrium with two non-zero coordinates, $z_{i-1}$ and $z_i$, to an equilibrium with only one non-zero coordinate~$z_{i+1}$.
For this reason, we label coordinates starting at $i-1$ rather than~$i$.
Figure~\ref{fig:plane-to-axis} illustrates important points in the construction of the global map.

We define the cross sections as follows:
 $$
 H_i^{out,i+1}=\{ (z_{i-1},z_i,z_{i+1},z_{i+2}) \in \R^4: \;\; z_{i-1}=\Order(1), z_i=\Order(1), z_{i+1}=h, z_{i+2}<h \}
 $$
and
 $$
 H_{i+1}^{in,i}=\{ (z_{i-1}, z_i,z_{i+1},z_{i+2}) \in \R^4: \;\; \max(z_{i-1},z_i)=h, z_{i+1}=\Order(1), z_{i+2}<h \},
 $$
for small $h$. 
At~$\xi_i$, $z_{i-1}$ and $z_i$ are the radial coordinates (and are of order~$1$) and $z_{i+1}$ is the expanding coordinate, and at~$\xi_{i+1}$, $z_{i+1}$~is the radial coordinate, $z_{i-1}$ and $z_i$ are the contracting coordinates, and $z_{i+2}$~is the expanding coordinate.
The reason for writing $\max(z_{i-1},z_i)=h$ is that the two contracting coordinates decay at different rates, so we define the cross section in terms of the one that reaches~$h$ last.
This guarantees that both $z_{i-1}\leq h$ and $z_i\leq h$ on the incoming section.
See the discussion around Figure~\ref{fig:logh} for more detail.

Again the unstable manifold of $\xi_i$, $W^u(\xi_i)$, is the connection from $\xi_i$ to $\xi_{i+1}$, and intersects the cross sections at
$$
W^u(\xi_i) \cap H_i^{out,i+1} = \{ (\bar{z}_{i-1},\bar{z}_i,h,0) \}
$$
and either
$$
W^u(\xi_i) \cap H_{i+1}^{in,i} = \{ (h,\hat{z}_i,\hat{z}_{i+1},0) \}
\qquad\mbox{or}\qquad 
W^u(\xi_i) \cap H_{i+1}^{in,i} = \{ (\hat{z}_{i-1},h,\hat{z}_{i+1},0) \},
$$
depending on which contracting eigenvalue is closer to zero.
The $P_i$ space ($z_{i+2}=0$) is invariant so $W^u(\xi_i)$ leaves $\xi_i$ with $z_{i+2}=0$ and arrives at $\xi_{i+1}$ also with $z_{i+2}=0$.
We write a point in $H_i^{out,i+1}$ as $(z_{i-1},z_i,z_{i+1},z_{i+2})=(\bar{z}_{i-1}+\tilde{z}_{i-1},\bar{z}_i+\tilde{z}_i,h,z_{i+2})$, with $|\tilde{z}_{i-1}|<h$ and $|\tilde{z}_i|<h$.
The linearisation around $W^u(\xi_i)$ provides the global map:
\begin{equation} \label{eq:globalplanetoaxis}
    \Psi_{i \to i+1} \left(\begin{array}{c} 
 z_{i-1} \\ z_i\\ z_{i+2}
    \end{array} \right)
    = \left(\begin{array}{ccc} 
[...] & [...] & [...] \\
\mbox{}[...] & [...] & [...] \\
0 & 0 & A^{i+2}_{i \to i+1}
    \end{array} \right) \left(\begin{array}{c} 
\tilde{z}_{i-1}\\ \tilde{z}_{i}\\ z_{i+2}
    \end{array} \right) + \left(\begin{array}{c} 
\hat{z}_{i}\\ \hat{z}_{i+1}\\ 0
    \end{array} \right),
\end{equation}
in the case that $z_{i-1}$ reaches~$h$ last.
This map takes as argument the values of $(z_{i-1},z_i,z_{i+2})$ on the outgoing section and returns the values of
$(z_i,z_{i+1},z_{i+2})$ on the incoming section.
In the case that $z_i$ reaches~$h$ last, the returned values would be $(z_{i-1},z_{i+1},z_{i+2})$ on the incoming section, with a $\hat{z}_{i-1}$ (instead of~$\hat{z}_{i}$) in the first line.
The invariance of the $z_{i+2}=0$ subspace leads to the structure of the $i+2$ part of the matrix.
The $A$~coefficient comes from the linearisation around the unstable manifold.
The un-named entries in the matrix, indicated by~$[...]$, lead to contributions that are small compared to the fixed order~$1$ values of~$\hat{z}_i$ and~$\hat{z}_{i+1}$.
There are also $\Order(h^2)$ corrections to this linearised map (not written).
Both of these small contributions to the global map will be disregarded when all variables are small, close to the heteroclinic cycle.

\begin{figure}
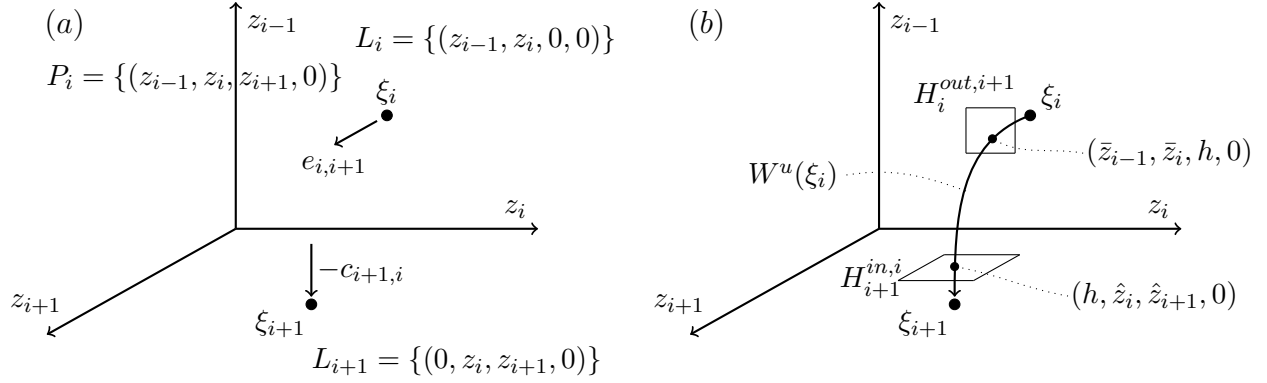

\hbox to \hsize{\fbox{%
\beginpgfgraphicnamed{Fig-global-map_PP_1}%
\endpgfgraphicnamed}%
\hfil\fbox{%
\beginpgfgraphicnamed{Fig-global-map_PP_2}%
\endpgfgraphicnamed}}
\caption{The global map between two equilibria in planes. 
Panel~(a) shows the flow-invariant subspaces, equilibria, and the direction of the eigenvalues. 
Panel~(b) shows the cross sections and coordinates of their intersection with the unstable manifold of~$\xi_i$.
\label{fig:plane-to-plane}}
\end{figure}

\paragraph{Plane-to-plane connection.}
This is the type of global map that describes the dynamics along a connection between two equilibria in different (but intersecting) coordinate planes. 
Such a connection appears, for example, between $\xi_2$ and $\xi_3$ in Case~3, Figure~\ref{fig:case3}(a). 
In this case, each equilibrium has two non-zero coordinates: $\xi_i$ has nonzero $z_{i-1}$ and~$z_i$, and $\xi_{i+1}$ has nonzero $z_i$ and~$z_{i+1}$.
As in the plane-to-axis case, we label coordinates starting at $i-1$ rather than~$i$.
Figure~\ref{fig:plane-to-plane} illustrates important points in the construction of the global map.

We define the cross sections as follows:
 $$
 H_i^{out,i+1}=\{ (z_{i-1},z_i,z_{i+1},z_{i+2}) \in \R^4: \;\; z_{i-1}=\Order(1), z_i=\Order(1), z_{i+1}=h, z_{i+2}<h \}
 $$
and
 $$
 H_{i+1}^{in,i}=\{ (z_{i-1},z_i,z_{i+1},z_{i+2}) \in \R^4: \;\; z_{i-1}=h, z_{i}=\Order(1), z_{i+1}=\Order(1), z_{i+2}<h \},
 $$
for small $h$. 
At~$\xi_i$, $z_{i-1}$ and $z_i$ are the radial coordinates (and are of order~$1$) and $z_{i+1}$ is the expanding coordinate, and at $\xi_{i+1}$, $z_i$ and $z_{i+1}$ are the radial coordinates, $z_{i-1}$ is the contracting
coordinate, and $z_{i+2}$~is the expanding coordinate.

Again the unstable manifold of $\xi_i$, $W^u(\xi_i)$, is the connection from $\xi_i$ to $\xi_{i+1}$, and intersects the cross sections at
$$
W^u(\xi_i) \cap H_i^{out,i+1} = \{ (\bar{z}_{i-1},\bar{z}_i,h,0) \}
\qquad\mbox{and}\qquad 
W^u(\xi_i) \cap H_{i+1}^{in,i} = \{ (h,\hat{z}_i,\hat{z}_{i+1},0) \}.
$$
The $P_i$ space ($z_{i+2}=0$) is invariant so $W^u(\xi_i)$ leaves $\xi_i$ with $z_{i+2}=0$ and arrives at $\xi_{i+1}$ also with $z_{i+2}=0$.
We write a point in $H_i^{out,i+1}$ as $(z_{i-1},z_i,z_{i+1},z_{i+2})=(\bar{z}_{i-1}+\tilde{z}_{i-1}, \bar{z}_i+\tilde{z}_i,h,z_{i+2})$, with $|\tilde{z}_{i-1}|<h$ and $|\tilde{z}_i|<h$.
The linearisation around $W^u(\xi_i)$ provides the global map
\begin{equation} \label{eq:globalplanetoplane}
    \Psi_{i \to i+1} \left(\begin{array}{c} 
 z_{i-1} \\ z_i\\ z_{i+2}
    \end{array} \right)
    = \left(\begin{array}{ccc} 
[...] & [...] & [...] \\
\mbox{}[...] & [...] & [...] \\
0 & 0 & A^{i+2}_{i \to i+1}
    \end{array} \right) \left(\begin{array}{c} 
\tilde{z}_{i-1}\\ \tilde{z}_{i}\\ z_{i+2}
    \end{array} \right) + \left(\begin{array}{c} 
\hat{z}_{i}\\ \hat{z}_{i+1}\\ 0
    \end{array} \right).
\end{equation}
This map takes as argument the values of $(z_{i-1},z_i,z_{i+2})$ on the outgoing section and returns the values of $(z_i,z_{i+1},z_{i+2})$ on the incoming section.                                          
The invariance of the $z_{i+2}=0$ subspace leads to the structure of the $i+2$ part of the matrix.                         The $A$~coefficient comes from the linearisation around the unstable manifold.
The un-named entries in the matrix, indicated by~$[...]$, lead to contributions that are small compared to the fixed order~$1$ values of~$\hat{z}_i$ and~$\hat{z}_{i+1}$.            
There are also $\Order(h^2)$ corrections to this linearised map (not written).
Both of these small contributions to the global map will be disregarded when all variables are small, close to the heteroclinic cycle.                                     

\smallskip

These global maps are all characterised as having two different parts, corresponding to the different natures of the upper and lower rows in the matrices.
The upper part is one row in connections starting on an axis and two rows in connections starting on a plane.
The upper parts are dominated by fixed $\Order(1)$ numbers (such as $\hat{z}_{i-1}$ in~\eqref{eq:globalaxistoaxis}), with $\Order(h)$ contributions, indicated by~$[...]$.
The exception is the second contracting direction in~\eqref{eq:globalplanetoaxis}, in which $\hat{z}_i<h$ and all the other terms are $\Order(h)$ or smaller.
We will see in the next subsections that all these terms can be neglected when considering properties of trajectories very close to the heteroclinic cycle, since $h$ is a \emph{fixed} small number, while trajectories can come arbitrarily close the cycle.

The lower parts of the global maps are all written in terms of $2\times2$ (starting on an axis) or $1\times1$ (starting on a plane) diagonal matrices. 
With~$\xi_i$ on an axis, the global map starting at $H_i^{out,i+1}$ is of the form
    $$
    \left( \begin{array}{c}
    z_{i+2} \\
    z_{i+3}
    \end{array} \right)
    \rightarrow
    \left( \begin{array}{c}
    A^{i+2}_{i \to i+1} z_{i+2} \\
    A^{i+3}_{i \to i+1} z_{i+3}
    \end{array} \right) + 
    \left( \begin{array}{c}
    \Order(h^2) \\
    \Order(h^2)
    \end{array} \right).
    $$
With $\xi_i$ on a plane, the global map starting at $H_i^{out,i+1}$ is of the form
    $$
    \left( \begin{array}{c}
    z_{i+2}
    \end{array} \right)
    \rightarrow
    \left( \begin{array}{c}
    A^{i+2}_{i \to i+1} z_{i+2}
    \end{array} \right) +
    \left( \begin{array}{c}
    \Order(h^2)
    \end{array} \right).
    $$
We have explicitly written the size of the nonlinear corrections to the maps.
However, note that the correction to (for example) $z_{i+2}$ must be zero when $z_{i+2}=0$ (because of the invariance of the $z_{i+2}=0$ space), so the correction can be thought of as being $\Order(hz_{i+2})$, and similarly for $z_{i+3}$.
These amount to $\Order(h)$ corrections to the $A$ coefficients for trajectories very close to the heteroclinic cycle.

\subsection{Local maps}\label{sec:localmaps}

The local map at equilibrium point~$\xi_i$, denoted by~$\phi_i$, takes the trajectory from $H_i^{in,i-1}$ to $H_i^{out,i+1}$, and describes the dynamics near~$\xi_i$.
As in~\cite{Podvigina2012}, we use logarithms of the coordinates, so that the local maps take the form of a linear map.
We assume (as is usual) that the flow near each equilibrium point is linearisable and depends only on the four eigenvalues of the Jacobian matrix at that point.

There are two types of local maps, depending on whether the equilibrium is on an axis or in a plane.
On an axis, there is one negative radial eigenvalue, one positive expanding eigenvalue, and two negative eigenvalues that can be either contracting or transverse.
In a plane, there are two negative radial eigenvalues, one positive expanding eigenvalue, and one negative eigenvalue that can be either contracting or transverse.
Whether an eigenvalue is contracting or transverse depends on the preceding global dynamics.

We proceed in the usual manner.
Near each equilibrium point $\xi_i$ we use the following linear approximation for the local expanding dynamics
 \begin{equation*}
 \dot{z}_{i+1} = e_{i,i+1}z_{i+1},
 \end{equation*}
where $e_{i,i+1}$ is the expanding eigenvalue at~$\xi_i$ in the $z_{i+1}$~direction.
The solution of this equation is $z_{i+1}(t)=z_{i+1}(0)\exp(e_{i,i+1}t)$, where $z_{i+1}(0)$ is the value of the expanding coordinate on~$H_{i}^{in,i-1}$.
The trajectory reaches~$H_{i}^{out,i+1}$ at time~$T$:
 \begin{equation}\label{eq:time}
 T = -\frac{1}{e_{i,i+1}}\log\left(\frac{z_{i+1}(0)}{h}\right),
 \end{equation}
where this is found from solving $z_{i+1}(T)=h$.

The radial, contracting and transverse directions all have negative eigenvalues, from~(C1).
We use~$z(t)$ to represent any of these, with eigenvalue~$-k$, so $-k$ is a radial, contracting or transverse eigenvalue.
The differential equation $\dot{z}=-kz$ has solution $z(t)=z(0)\exp(-kt)$.
In the contracting and transverse directions, the invariance of the subspaces means that eigenvalues must be real, by~(A2).

In the radial case, $z$~represents the deviation from the equilibrium point, and when there are two radial eigenvalues, these can be complex.
However, as in~\cite{Krupa1995}, the radial eigenvalues are irrelevant to the stability of the heteroclinic cycle, even in the absence of a contracting eigenvalue.
The reason is that we only need to know that the radial coordinate is~$\Order(1)$, not its exact value, since radial coordinates at one point become contracting coordinates at the next.
Incoming cross-sections are defined by requiring that the contracting coordinate (if there is only one) is equal to~$h$.
If there is more than one, we show below that they are of similar size (defined more precisely below).
The end result is that we do not need the exact form of the local maps in the radial direction.

The value of~$z$ at time~$T$, when the trajectory reaches~$H_{i}^{out,i+1}$, is
 $$
 z(T) = z(0)\left(\frac{z_{i+1}(0)}{h}\right)^{k/e_{i,i+1}}.
 $$
Writing this in terms of logarithms, we have
 \begin{equation}\label{eq:generallogmap}
 \log z(T) = \log z(0) + \frac{k}{e_{i,i+1}}\log z_{i+1}(0) - \frac{k}{e_{i,i+1}}\log h.
 \end{equation}
Recall that $h$ is a fixed small number, but near the heteroclinic cycle, $z_{i+1}(0)$ is arbitrarily small,
so $|\log z_{i+1}(0)| \gg |\log h|$, taking absolute values as both logarithms are negative.
As a result, we write the last term in~\eqref{eq:generallogmap} as~$\Order(\log h)$.
This term is unimportant when trajectories are very close to the cycle.

We next turn to the $\log z(0)$ term in~\eqref{eq:generallogmap}: we treat this term differently according to whether $z$~is a transverse or contracting coordinate.
If $z$ is a transverse coordinate, we have $z(0)\ll h$ and so $|\log z(0)|\gg|\log h|$, and the $\log z(0)$ term must be retained. 
Conversely, when $z$ is a contracting direction, this coordinate was $\Order(1)$ at some point along the connecting trajectory prior to reaching~$H_i^{in,i-1}$.
If there is only one contracting direction, in which case $z(0)=h$, the first term on the RHS of~\eqref{eq:generallogmap} is~$\log h$ and is also unimportant when trajectories are very close to the cycle.
However, in cycles in pluridimensions, there can be equilibria with more than one contracting direction.
In this case, we argue in the next paragraph that $|\log z(0)|=\Order(|\log h|)$, and so the $\log z(0)$ term can be absorbed into the $\Order(\log h)$ term in~\eqref{eq:generallogmap}, and so the contracting directions are all treated in the same way.

\begin{figure}
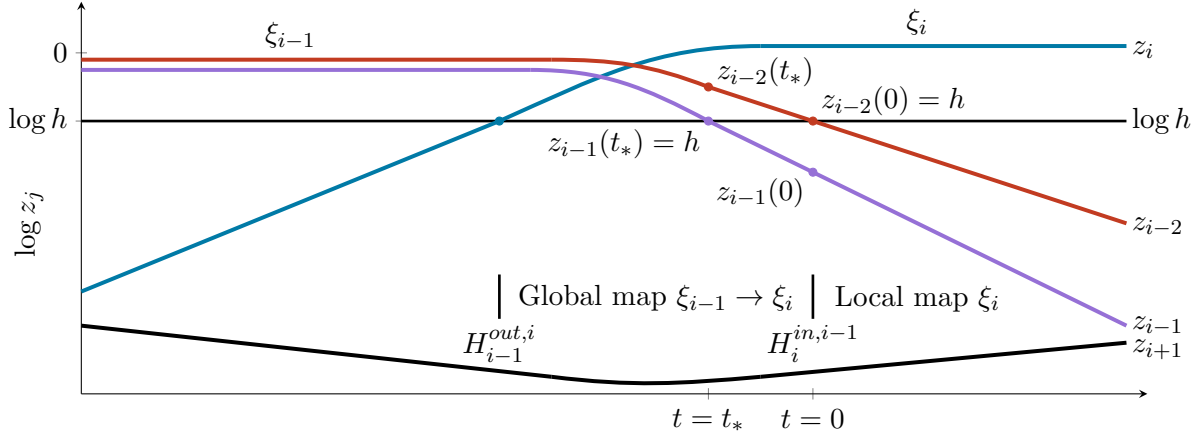
  
\hbox to \textwidth{\hfil\fbox{%
\beginpgfgraphicnamed{x1x2x3x4_timeseries_log_h}%
\endpgfgraphicnamed}\hfil}%
\caption{Illustration of the sizes of the variables as a function of time at~$H_i^{in,i-1}$, 
the incoming section of the local map at~$\xi_i$.
The location of the cross-sections is indicated above the horizontal axis with a reference to the global and local maps, ending or starting at~$H_i^{in,i-1}$, respectively.
The colors represent $z_{i-2}$~(red), $z_{i-1}$~(purple), $z_i$~(blue) and $z_{i+1}$~(black), plotted in logarithmic coordinates. 
The local map begins at time $t=0$, with $\max(z_{i-2}(0),z_{i-1}(0))=h$.
In this illustration, we take $z_{i-2}(0)=h$ so $z_{i-1}(0)<h$,
and there is an earlier time $t_*<0$ such that $z_{i-1}(t_*)=h$ and $z_{i-2}(t_*)>h$.
This happens because the previous equilibrium point~$\xi_{i-1}$ had both $z_{i-2}$ and $z_{i-1}$
being of order one. 
Throughout, the fourth variable~$z_{i+1}$ is small compared to~$h$.
\label{fig:logh}}
\end{figure}

In the case where there is more than one contracting direction, the trajectory intersects the incoming cross section when the largest of the contracting coordinates is equal to~$h$ (i.e., $\max(z_{i-2}(0),z_{i-1}(0))=h$ in analogy with the plane-to-axis discussion in Section~\ref{sec:globalmaps} above).
To be definite, as in that case, we suppose that $z_{i-2}$ reaches~$h$ last, so $z_{i-2}(0)=h$ and $z_{i-1}(0)<h$.
We define an earlier time $t_*<0$ such that at that time, $z_{i-1}(t_*)=h$ and $z_{i-2}(t_*)>h$.
See Figure~\ref{fig:logh} for more detail.
We take $h$ small enough so that the dynamics is governed by linear differential equations from time~$t_*$, and so
 \[
 z_{i-2}(t)=h e^{-c_{i-2}(t-0)}
 \qquad\text{and}\qquad
 z_{i-1}(t)=h e^{-c_{i-1}(t-t_*)},
 \]
where $-c_{i-2}$ and $-c_{i-1}$ are the relevant contracting eigenvalues.
This gives
 \[
 z_{i-2}(t_*)=h e^{-c_{i-2}t_*}
 \qquad\text{and}\qquad
 z_{i-1}(0)=h e^{c_{i-1} t_*}.
 \]
We eliminate $t_*$ between these expressions and use the fact that $0>\log z_{i-2}(t_*)>\log h$ to derive
\[
\log h > \log z_{i-1}(0) > \left(1+\frac{c_{i-1}}{c_{i-2}}\right)\log h.
\]
From this we conclude that $\log z_{i-1}(0)=\Order(\log h)$ for any contracting coordinate, and so the first term in~\eqref{eq:generallogmap} can be absorbed into the $\Order(\log h)$ terms.

This distinction between contracting and transverse coordinates arises because, in the contracting case, the coordinates have $\Order(1)$ values at the previous equilibrium point, while transverse coordinates are generally very small.
This distinction becomes important when we compose the local and global maps in Section~\ref{sec:transitionmatrices}.

\paragraph{Equilibrium on an axis.}
Here we consider the local map~$\phi_i$ from $H_{i}^{in,i-1}$ to $H_{i}^{out,i+1}$, at an equilibrium point $\xi_i$ on an axis.
At~$\xi_i$, the four directions are radial~($z_i$), expanding~($z_{i+1}$), contracting~($z_{i-1}$) and a fourth direction ($z_{i-2}$).
This last direction can be contracting or transverse, so we label its eigenvalue as~$(ct)_{i,i-2}$.
In terms of logarithms, we have
    \begin{equation} \label{eq:maplocalaxis}
    \phi_i
    \left(\begin{array}{c}
    \log z_{i+1}\\
    \log z_{i-2} 
    \end{array} \right)
    =
    \left[ \begin{array}{cc}
    \frac{\displaystyle(ct)_{i,i-2}}{\displaystyle e_{i,i+1}} & 1 \\[8pt]
    \frac{\displaystyle c_{i,i-1}}{\displaystyle e_{i,i+1}} & 0
    \end{array} \right]
    \left(\begin{array}{c}
    \log z_{i+1}\\
    \log z_{i-2} 
    \end{array} \right)
     + \Order(\log h).
    \end{equation}
This map takes as argument the logarithms of $(z_{i+1},z_{i-2})$ on the incoming section and returns the logarithms of $(z_{i-2},z_{i-1})$ on the outgoing section.
In the case where $z_{i-2}$ is a contracting coordinate, the discussion of having more than one contracting coordinate applies, and so $\log z_{i-2}=\Order(\log h)$, and the terms that arise from the second column in the matrix are~$\Order(\log h)$.

\paragraph{Equilibrium in a plane.}
Here we consider the local map~$\phi_i$ from $H_{i}^{in,i-1}$ to $H_{i}^{out,i+1}$, at an equilibrium point $\xi_i$ in a plane.
At~$\xi_i$, the four directions are radial~($z_{i-1}$ and $z_i$), expanding~($z_{i+1}$) and one more direction ($z_{i-2}$).
This last direction can be contracting or transverse, so we label its eigenvalue as~$(ct)_{i,i-2}$.
In terms of logarithms, we have
    \begin{equation} \label{eq:maplocalplane}
    \phi_i
    \left(\begin{array}{c}
    \log z_{i+1}\\
    \log z_{i-2} 
    \end{array} \right)
    =
    \left[ \begin{array}{cc}
    \frac{\displaystyle(ct)_{i,i-2}}{\displaystyle e_{i,i+1}} & 1 
    \end{array} \right]
    \left(\begin{array}{c}
    \log z_{i+1}\\
    \log z_{i-2} 
    \end{array} \right)
     + \Order(\log h).
    \end{equation}
This map takes as argument the logarithms of $(z_{i+1},z_{i-2})$ on the incoming section and returns the logarithm of $z_{i-2}$ on the outgoing section.
As above, in the case where $z_{i-2}$ is a contracting coordinate, the term that arises from the second column in the matrix is~$\Order(\log h)$.

\subsection{Transition matrices: local maps composed with global maps}\label{sec:transitionmatrices}

We next compose local and global maps, $\Psi_{i\to i+1}\circ\phi_i$, going from $H_i^{in,i-1}$ to $H_{i+1}^{in,i}$.
The local maps are written in terms of the logarithms of the coordinates, and we use the same logarithmic representation for the composed maps.
The local map from $H_i^{in,i-1}$ to $H_{i}^{out,i+1}$ depends not only on the position (axis or plane) of $\xi_i$ but also on the position of~$\xi_{i-1}$. 
The reason for this dependence is that, at $H_i^{in,i-1}$, $z_{i-2}$ can be a contracting or transverse coordinate, according to whether $\xi_{i-1}$ is in a plane or on an axis.
This has a consequence for whether $\log z_{i-2}=\Order(\log h)$ (and so can be neglected) or $|\log z_{i-2}|\gg|\log h|$ (and so must be kept).
There are thus eight cases for $\xi_{i-1}\rightarrow\xi_i\rightarrow\xi_{i+1}$: axis-to-axis-to-axis, etc.
However, the form of the global map from $\xi_i\rightarrow\xi_{i+1}$ depends on the position of $\xi_i$ (and not~$\xi_{i+1}$), so there are in fact only four distinct cases for the composed local and global maps.

\begin{figure}
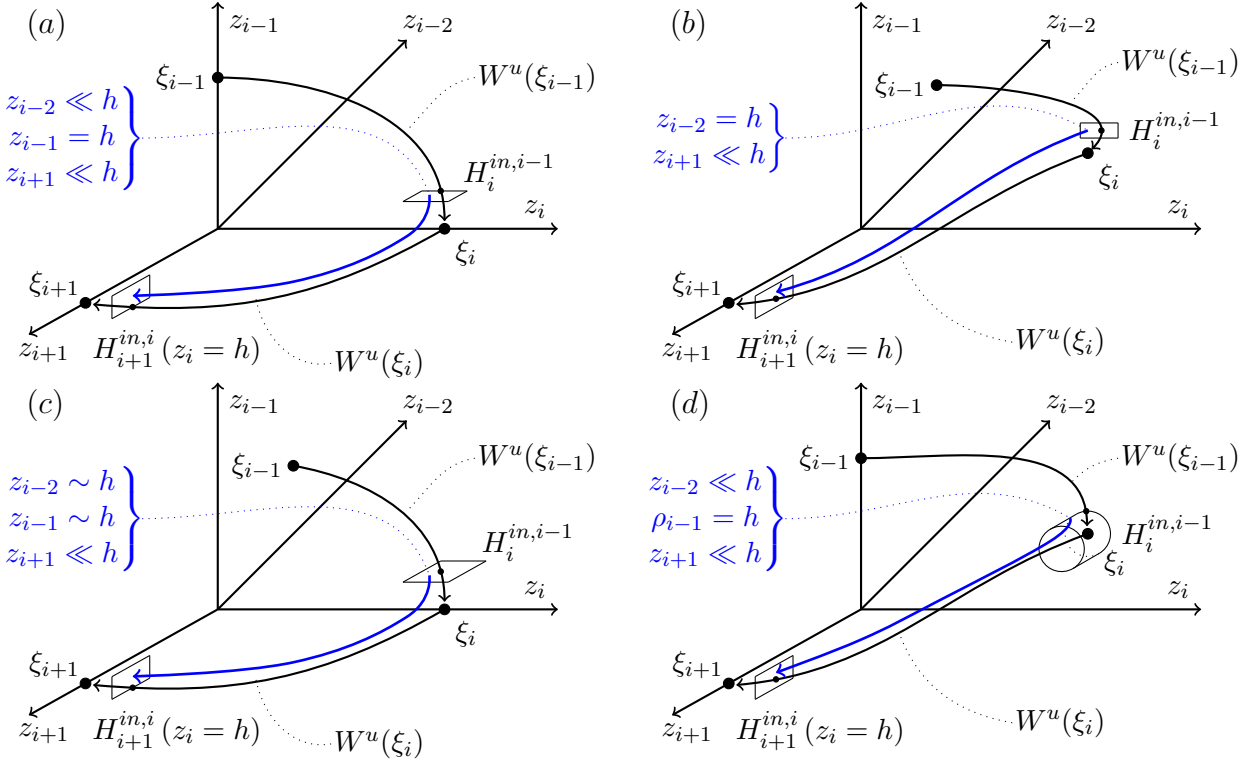

\hbox to \hsize{\fbox{%
                \beginpgfgraphicnamed{Fig-composed-maps_AAA}%
                \endpgfgraphicnamed}\hfil%
                \fbox{%
                \beginpgfgraphicnamed{Fig-composed-maps_PPA}%
                \endpgfgraphicnamed}}
\hbox to \hsize{\fbox{%
                \beginpgfgraphicnamed{Fig-composed-maps_PAA}%
                \endpgfgraphicnamed}\hfil%
                \fbox{%
                \beginpgfgraphicnamed{Fig-composed-maps_APA}%
                \endpgfgraphicnamed}}
\caption{Four examples of composed maps: (a)~axis-to-axis-to-axis, (b)~plane-to-plane-to-axis, (c)~plane-to-axis-to-axis and (d)~axis-to-plane-to-axis.
The maps are all from $H_{i}^{in,i-1}$ to $H_{i+1}^{in,i}$, trajectories having come from~$\xi_{i-1}$.
We have presented all four dimensions, with $z_{i+1}$~coming out of the page and $z_{i-2}$~going into the page.
The black curves indicate $W^u(\xi_{i-1})$ and~$W^u(\xi_i)$.
The blue curve indicates a trajectory from $H_{i}^{in,i-1}$ to $H_{i+1}^{in,i}$.
The composed map is given by \eqref{eq:mapAAX} in panel~(a), \eqref{eq:mapPPX}~in panel~(b), \eqref{eq:mapPAX}~in panel~(c), and \eqref{eq:mapAPX}~in panel~(d).
The magnitudes of the variables at $H_{i}^{in,i-1}$ are indicated: $z\sim h$ means $\log z=\Order(\log h)$ (it could also mean $z=h$).
Variables that are not given are~$\Order(1)$.
\label{fig:XtoYtoZ}}
\end{figure}

The composed maps (see Figure~\ref{fig:XtoYtoZ}) will include $\Order(\log h)$ and other terms that can be neglected for trajectories that are very close to the heteroclinic cycle.
Once these terms are neglected, the result is a transition matrix~\cite{Field1991,Krupa2004} that describes the map from $H_i^{in,i-1}$ to~$H_{i+1}^{in,i}$.
The transition matrix multiplies the logarithms of the small (expanding and transverse) coordinates on~$H_i^{in,i-1}$.
Note that our transition matrices are not necessarily square: this is a new feature of cycles in pluridimensions, explained in more detail below.

We use the notation of~\cite{Podvigina2011} and define $a_i^{(j)}$ to be the negative of the quotient between the contracting (in the $j$~direction) and expanding eigenvalues at~$\xi_i$.
We define $b_i^{(j)}$ to be the negative of the quotient between the transverse (in the $j$~direction) and expanding eigenvalues at~$\xi_i$. 

\paragraph{Axis-to-axis-to-axis and axis-to-axis-to-plane cases.} 
Here we consider the local map from $H_{i}^{in,i-1}$ to $H_{i}^{out,i+1}$, at an equilibrium point $\xi_i$ on an axis, having come from $\xi_{i-1}$, also on an axis.
The next point~$\xi_{i+1}$ can be on an axis or on a plane.
This situation is illustrated in Case~2, Figure~\ref{fig:case2}: 
$\xi_3\rightarrow\xi_4\rightarrow\xi_1$ is axis-to-axis-to-axis, and
$\xi_4\rightarrow\xi_1\rightarrow\xi_2$ is axis-to-axis-to-plane.

At~$\xi_i$, the four directions are radial~($z_i$), contracting~($z_{i-1}$), expanding~($z_{i+1}$) and transverse ($z_{i-2}$), with $z_i=\Order(1)$ and $z_{i-1}=h$, and $z_{i+1}\ll h$ and $z_{i-2}\ll h$.
After applying the local map $\phi_i$~\eqref{eq:maplocalaxis}, we get
 $$
    \left[ \begin{array}{cc}
    \frac{\displaystyle t_{i,i-2}}{\displaystyle e_{i,i+1}} & 1 \\[8pt]
    \frac{\displaystyle c_{i,i-1}}{\displaystyle e_{i,i+1}} & 0
    \end{array} \right]
    \left(\begin{array}{c}
    \log z_{i+1}\\
    \log z_{i-2} 
    \end{array} \right)
     + \Order(\log h).
 $$
This gives, in the first row, $\log z_{i-2}$ on $H_{i}^{out,i+1}$, and in the second row, $\log z_{i-1}$.

The next step is to go from $\xi_i$ (on an axis) to $\xi_{i+1}$, which could be on an axis or in a plane.
The global map from $H_{i}^{out,i+1}$ to $H_{i+1}^{in,i}$, in both the axis to axis case~\eqref{eq:globalaxistoaxis} and in the axis to plane case~\eqref{eq:globalaxistoplane}, is written in terms of $z_{i+2}$ and~$z_{i+3}$.
Since we are working in $\R^4$, we have $z_{i+2}=z_{i-2}$ and~$z_{i+3}=z_{i-1}$, but in higher dimensions, the maps would have to keep track of additional variables.
In~\eqref{eq:globalaxistoaxis} and ~\eqref{eq:globalaxistoplane}, both maps act to multiply $z_{i+2}$ and~$z_{i+3}$ by order~$1$ constants $A_{i\rightarrow i+1}^{i+2}$ and~$A_{i\rightarrow i+1}^{i+3}$.
Composing the local and global maps results in
\begin{equation} \label{eq:mapAAX}
    \left[ \begin{array}{cc}
    b_i^{(i-2)} & 1 \\[8pt]
    a_i^{(i-1)} & 0
    \end{array} \right]
    \left(\begin{array}{c}
    \log z_{i+1}\\
    \log z_{i-2} 
    \end{array} \right)
     + \Order(\log h)
     + \Order(\log A),
\end{equation}
where $A$ represents $A_{i\rightarrow i+1}^{i+2}$ and~$A_{i\rightarrow i+1}^{i+3}$.
This gives, in the first row, $\log z_{i+2}$ on $H_{i+1}^{in,i}$, and in the second row, $\log z_{i+3}$.
Once the $\Order(\log h)$ and $\Order(\log A)$ terms are neglected, the map~\eqref{eq:mapAAX} gives the transition matrix for the axis-to-axis-to-axis and axis-to-axis-to-plane cases.

\paragraph{Plane-to-plane-to-axis and plane-to-plane-to-plane cases.}
Here we consider the local map from $H_{i}^{in,i-1}$ to $H_{i}^{out,i+1}$, at an equilibrium point $\xi_i$ on an plane, having come from $\xi_{i-1}$, also on an plane.
The next point~$\xi_{i+1}$ can be on an axis or on a plane.
This situation is illustrated in Case~3, Figure~\ref{fig:case3}: $\xi_2\to\xi_3\to\xi_4$ is plane-to-plane-to-axis, and in Case~4, Figure~\ref{fig:case4}: $\xi_2\to\xi_3\to\xi_4$ is plane-to-plane-to-plane.

At $\xi_i$, the four directions are radial~($z_i$, $z_{i-1}$), contracting~($z_{i-2}$) and expanding~($z_{i+1}$), with no transverse direction, with $z_{i-1},z_i=\Order(1)$ and $z_{i-2}=h$, and $z_{i+1}\ll h$.
After applying the local map $\phi_i$~\eqref{eq:maplocalplane}, we get
 $$
    \left[ \begin{array}{cc}
    \frac{\displaystyle c_{i,i-2}}{\displaystyle e_{i,i+1}} & 1 
    \end{array} \right]
    \left(\begin{array}{c}
    \log z_{i+1}\\
    \log z_{i-2} 
    \end{array} \right)
     + \Order(\log h).
$$
Since $z_{i-2}=h$ on~$H_{i}^{in,i-1}$, the $\log z_{i-2}$ contribution is absorbed into the $\Order(\log h)$ term, and so the right column of the matrix is removed.
This gives $\log z_{i-2}$ on $H_{i}^{out,i+1}$.

The next step is to go from $\xi_i$ (on an plane) to $\xi_{i+1}$, which could be on an axis or in a plane.
The global map from $H_{i}^{out,i+1}$ to $H_{i+1}^{in,i}$, in both the plane to axis case~\eqref{eq:globalplanetoaxis} and in the plane to plane case~\eqref{eq:globalplanetoplane}, is written in terms of~$z_{i+2}$.
Since we are working in $\R^4$, we have $z_{i+2}=z_{i-2}$.           
In~\eqref{eq:globalplanetoaxis} and~\eqref{eq:globalplanetoplane}, both maps act to multiply $z_{i+2}$ by an order~$1$ constant~$A_{i\rightarrow i+1}^{i+2}$.
Composing the local and global maps results in
    \begin{equation} \label{eq:mapPPX}
    \left[ \begin{array}{c}
    a_i^{(i-2)}
    \end{array} \right]
    \left(\begin{array}{c}
    \log z_{i+1} 
    \end{array} \right)
     + \Order(\log h) + \Order(\log A),
    \end{equation}
where $A$ represents~$A_{i\rightarrow i+1}^{i+2}$.
The first (only) row is $z_{i+2}$ on $H_{i+1}^{in,i}$.
Once the $\Order(\log h)$ and $\Order(\log A)$ terms are neglected, the map~\eqref{eq:mapPPX} gives the transition matrix for the plane-to-plane-to-axis and plane-to-plane-to-plane cases.

\paragraph{Plane-to-axis-to-axis and plane-to-axis-to-plane cases.}
Here we consider the local map from $H_{i}^{in,i-1}$ to $H_{i}^{out,i+1}$, at an equilibrium point $\xi_i$ on an axis, having come from $\xi_{i-1}$ on a plane.
The next point~$\xi_{i+1}$ can be on an axis or on a plane.
This situation is illustrated in Case~2, Figure~\ref{fig:case2}: 
$\xi_2\rightarrow\xi_3\rightarrow\xi_4$ is plane-to-axis-to-axis, and
Case~1, Figure~\ref{fig:case1}: 
$\xi_2\rightarrow\xi_3\rightarrow\xi_4$ is plane-to-axis-to-plane.

At~$\xi_i$, the four directions are radial~($z_i$), contracting~($z_{i-1}$, $z_{i-2}$), expanding~($z_{i+1}$), with no transverse direction, with $z_i=\Order(1)$, $\max(z_{i-1},z_{i-2})=h$, and $z_{i+1}\ll h$.
After applying the local map $\phi_i$~\eqref{eq:maplocalaxis}, we get
 $$               
    \left[ \begin{array}{cc}
    \frac{\displaystyle c_{i,i-2}}{\displaystyle e_{i,i+1}} & 1 \\[8pt]
    \frac{\displaystyle c_{i,i-1}}{\displaystyle e_{i,i+1}} & 0
    \end{array} \right]
    \left(\begin{array}{c}
    \log z_{i+1}\\
    \log z_{i-2} 
    \end{array} \right)
     + \Order(\log h).
     $$
Since $\log z_{i-2}=\Order(\log h)$ on~$H_{i}^{in,i-1}$ as discussed in Section~\ref{sec:localmaps}, the $\log z_{i-2}$ contribution is absorbed into the $\Order(\log h)$ term, and so the right column of the matrix is removed.
This gives, in the first row, $\log z_{i-2}$ on $H_{i}^{out,i+1}$, and in the second row,~$\log z_{i-1}$.

The next step is to go from $\xi_i$ (on an axis) to $\xi_{i+1}$, which could be on an axis or in a plane.
The global map from $H_{i}^{out,i+1}$ to $H_{i+1}^{in,i}$, in both the axis to axis case~\eqref{eq:globalaxistoaxis} and in the axis to plane case~\eqref{eq:globalaxistoplane}, is written in terms of $z_{i+2}$ and~$z_{i+3}$.
Since we are working in $\R^4$, we have $z_{i+2}=z_{i-2}$ and~$z_{i+3}=z_{i-1}$.
In~\eqref{eq:globalaxistoaxis} and ~\eqref{eq:globalaxistoplane}, both maps act to multiply $z_{i+2}$ and~$z_{i+3}$ by order~$1$ constants $A_{i\rightarrow i+1}^{i+2}$ and~$A_{i\rightarrow i+1}^{i+3}$.
Composing the local and global maps results in
 \begin{equation} \label{eq:mapPAX}
    \left[ \begin{array}{c}
    a_i^{(i-2)} \\
    a_i^{(i-1)} 
    \end{array} \right]
    \left(\begin{array}{c}
    \log z_{i+1} 
    \end{array} \right)
     + \Order(\log h) + \Order(\log A),
    \end{equation}
where $A$ represents $A_{i\rightarrow i+1}^{i+2}$ and~$A_{i\rightarrow i+1}^{i+3}$.
This gives, in the first row, $\log z_{i+2}$ on $H_{i+1}^{in,i}$, and in the second row,~$\log z_{i+3}$.
Once the $\Order(\log h)$ and $\Order(\log A)$ terms are neglected, the map~\eqref{eq:mapPAX} gives the transition matrix for the plane-to-axis-to-axis and plane-to-axis-to-plane cases.

\paragraph{Axis-to-plane-to-axis and axis-to-plane-to-plane cases.} 
Here we consider the local map from $H_{i}^{in,i-1}$ to $H_{i}^{out,i+1}$, at an equilibrium point $\xi_i$ on an plane, having come from $\xi_{i-1}$ on an axis.
The next point~$\xi_{i+1}$ can be on an axis or on a plane.
This situation is illustrated in Case~1, Figure~\ref{fig:case1}: $\xi_1\to\xi_2\to\xi_3$ is axis-to-plane-to-axis, and in Case~3, Figure~\ref{fig:case3}: $\xi_1\to\xi_2\to\xi_3$ is axis-to-plane-to-plane.

At $\xi_i$, the four directions are radial~($z_i$, $z_{i-1}$), transverse~($z_{i-2}$) and expanding~($z_{i+1}$), with no contracting direction, with $z_{i-1},z_i=\Order(1)$ and $z_{i-2}\ll h$, and $z_{i+1}\ll h$.
In this case incoming section is a cylinder of radius~$h$.
After applying the local map $\phi_i$~\eqref{eq:maplocalplane}, we get
$$
    \left[ \begin{array}{cc}
    \frac{\displaystyle t_{i,i-2}}{\displaystyle e_{i,i+1}} & 1 
    \end{array} \right]
    \left(\begin{array}{c}
    \log z_{i+1}\\
    \log z_{i-2} 
    \end{array} \right)
     + \Order(\log h).
$$
This gives $\log z_{i-2}$ on $H_{i}^{out,i+1}$.

The next step is to go from $\xi_i$ (on an plane) to $\xi_{i+1}$, which could be on an axis or in a plane.
The global map from $H_{i}^{out,i+1}$ to $H_{i+1}^{in,i}$, in both the plane to axis case~\eqref{eq:globalplanetoaxis} and in the plane to plane case~\eqref{eq:globalplanetoplane}, is written in terms of~$z_{i+2}$.
Since we are working in $\R^4$, we have $z_{i+2}=z_{i-2}$.           
In~\eqref{eq:globalplanetoaxis} and~\eqref{eq:globalplanetoplane}, both maps act to multiply $z_{i+2}$ by an order~$1$ constant~$A_{i\rightarrow i+1}^{i+2}$.
Composing the local and global maps results in
\begin{equation}\label{eq:mapAPX}
    \left[ \begin{array}{cc}
    b_i^{(i-2)} & 1 
    \end{array} \right]
    \left(\begin{array}{c}
    \log z_{i+1}\\
    \log z_{i-2} 
    \end{array} \right)
     + \Order(\log h) + \Order(\log A),
\end{equation}
where $A$ represents~$A_{i\rightarrow i+1}^{i+2}$.
The first (only) row is $z_{i+2}$ on $H_{i+1}^{in,i}$.
Once the $\Order(\log h)$ and $\Order(\log A)$ terms are neglected, the map~\eqref{eq:mapAPX} gives the transition matrix for the axis-to-plane-to-axis and axis-to-plane-to-plane cases.

\subsection{Composing transition matrices around the cycle}
The four cases of cycles in pluridimensions in~$\R^4$ satisfying (A1)--(A4) are listed in Table~\ref{table:fourcases}.
Our choice of $\dim P_1=2$ and $\dim P_2=3$ means that $\xi_1$ is on an axis and $\xi_2$ is in a plane, so the local map at~$\xi_2$ is of the form~\eqref{eq:mapAPX}.
The implication is that on $H_2^{out,3}$, and on $H_3^{in,2}$, $z_4$~is the only small coordinate.
As a result, the Poincar\'e return map, composed around the whole cycle from $H_3^{in,2}$ to itself, takes the form of a $1\times1$ transition matrix in all four cases.
In the four cases, these matrices are:
 \begin{align}
 \label{eq:case1transition}
 &\text{Case 1:}&
    \left[ \begin{array}{c}
    \delta_1
    \end{array} \right]
 =&
    \left[ \begin{array}{cc}
    b_2^{(4)} & 1 
    \end{array} \right]
    \left[ \begin{array}{c}     
    a_1^{(3)} \\
    a_1^{(4)} 
    \end{array} \right]
    \left[ \begin{array}{cc}
    b_4^{(2)} & 1 
    \end{array} \right]
    \left[ \begin{array}{c}     
    a_3^{(1)} \\
    a_3^{(2)} 
    \end{array} \right],\\
 \label{eq:case2transition}
 &\text{Case 2:}&
    \left[ \begin{array}{c}
    \delta_2
    \end{array} \right]
 =&
    \left[ \begin{array}{cc}
    b_2^{(4)} & 1 
    \end{array} \right]
    \left[ \begin{array}{cc}
    b_1^{(3)} & 1 \\
    a_1^{(4)} & 0
    \end{array} \right]
    \left[ \begin{array}{cc}
    b_4^{(2)} & 1 \\
    a_4^{(3)} & 0
    \end{array} \right]
    \left[ \begin{array}{c}     
    a_3^{(1)} \\
    a_3^{(2)} 
    \end{array} \right],\\
 \label{eq:case3transition}
 &\text{Case 3:}&
    \left[ \begin{array}{c}
    \delta_3
    \end{array} \right]
 =&
    \left[ \begin{array}{cc}
    b_2^{(4)} & 1 
    \end{array} \right]
    \left[ \begin{array}{cc}
    b_1^{(3)} & 1 \\
    a_1^{(4)} & 0
    \end{array} \right]
    \left[ \begin{array}{c}     
    a_4^{(2)} \\
    a_4^{(3)} 
    \end{array} \right]
    \left[ \begin{array}{c}
    a_3^{(1)}
    \end{array} \right],\\
 \label{eq:case4transition}
 &\text{Case 4:}&
    \left[ \begin{array}{c}
    \delta_4
    \end{array} \right]
 =&
    \left[ \begin{array}{cc}
    b_2^{(4)} & 1 
    \end{array} \right]
    \left[ \begin{array}{c}     
    a_1^{(3)} \\
    a_1^{(4)} 
    \end{array} \right]
    \left[ \begin{array}{c}
    a_4^{(2)}
    \end{array} \right]
    \left[ \begin{array}{c}
    a_3^{(1)}
    \end{array} \right].
    \qquad\qquad\qquad\qquad
 \end{align}
Starting at a different cross section can result in $2\times2$ matrices, which have two eigenvalues: one the same as the~$\delta_i$ calculated starting at $H_3^{in,2}$ and the other equal to zero. 
For example, in Case~3, if we started at $H_2^{in,1}$ instead, the product of the matrices would be
 $$
    \left[ \begin{array}{cc}
    b_1^{(3)} & 1 \\
    a_1^{(4)} & 0
    \end{array} \right]
    \left[ \begin{array}{c}     
    a_4^{(2)} \\
    a_4^{(3)} 
    \end{array} \right]
    \left[ \begin{array}{c}
    a_3^{(1)}
    \end{array} \right]
    \left[ \begin{array}{cc}
    b_2^{(4)} & 1 
    \end{array} \right].
 $$
This is a $2\times2$ matrix with determinant equal to zero and with trace equal to~$\delta_3$, so the eigenvalues are zero and~$\delta_3$.
The $\delta_3$ eigenvalue determines the stability of the cycle, and retaining the $\Order(\log h)$ and $\Order(\log A)$ terms would break the degeneracy of the zero eigenvalue.

Each time around the cycle in Case~$i$, $\log z_4$ increases by a factor of~$\delta_i$, with corrections that are small compared with $|\log z_4|$ as $z_4$ goes to zero. 
Hence, the stability of the cycle in pluridimensions in Case~$i$ is given by~$\delta_i$: the cycle is asymptotically stable when $\delta_i>1$ and unstable when $\delta_i<1$ (recalling that all radial eigenvalues are negative), as is standard in stability calculations of heteroclinic cycles~\cite{Krupa1995}.

The time taken to go around the cycle is the sum of four short times jumping between the equilibria and four long times in the neighbourhoods of the four equilibria.
The total time is dominated by the four long times, which from~\eqref{eq:time} are proportional to the logarithms of the expanding coordinates at each point.
The logarithms of these expanding coordinates increase (in the stable case) by a factor of~$\delta_i$ each time around the cycle.
In the unstable case, the logarithms decrease by a factor of~$\delta_i$.
These variations in $\log z_4$ (in~$H_3^{in,2}$) and in the times taken to go around the cycle can be seen in the figures illustrating the examples in Section~\ref{sec:examples}.

\section{Specific examples of the four cases in $\R^4$}\label{sec:examples}

We construct ODEs for each of the four cases in Table~\ref{table:fourcases}, and illustrate the stability results by choosing two sets of parameter values in each case, with $\delta_i>1$ and $\delta_i<1$.
We use the same case labels as in Section~\ref{sec:pluri}.
The first of these examples (Case~1) is inspired by the convection and magnetoconvection examples of~\cite{Matthews1996,Rucklidge1995b}.
The stability of a similar example was considered by~\cite{Postlethwaite2005b}, with a similar calculation of the stability of an example of cycling chaos in~\cite{Ashwin1998d}.
The examples of the other cases are entirely new.

In this section, we use the specific coordinates $x_1$, \dots, $x_4$ rather than the general coordinates $z_{i-2}$, \dots, $z_{i+3}$.

\subsection{Case 1}

The first example is a cycle between the four equilibria $\xi_1=(1,0,0,0)$, $\xi_2=(d_1,1,0,0)$, $\xi_3=(0,0,1,0)$ and $\xi_4=(0,0,d_3,1)$, where we have scaled the four variables to set some components to be~$1$, and we leave the other components as parameters, with $d_1>0$ and $d_3>0$.
We set the eigenvalues at the origin to be $1$, $\pm 1$, $\pm 1$ and $\pm 1$ in the $x_1$, $x_2$, $x_3$ and $x_4$ directions respectively, scaling time so that the first eigenvalue is~$1$.
Requiring radial stability of the four equilibria leads us to the following choice of sign for the eigenvalues in the $x_2$, $x_3$ and $x_4$ directions: $-1$, $+1$ and $-1$ (these are chosen to be equal in magnitude to simplify the presentation).
Expressing the coefficients in the ODE in terms of the eigenvalues, we have
\begin{equation}\label{eq:ODE-example-1}
\begin{aligned}
\dot{x_1} & = x_1[\phantom{-}1-\mathcal{X}             + d_1x_2                     - c_{31}x_3 + (d_3(1+c_{31})+e_{41})x_4], \\
\dot{x_2} & = x_2[ -1+\mathcal{X} + e_{12}x_1 - d_1(1+e_{12})x_2           - c_{32}x_3 + (d_3(-1+c_{32})-t_{42})x_4],\\
\dot{x_3} & = x_3[\phantom{-}1-\mathcal{X} - c_{13}x_1 + (d_1(1+c_{13})+e_{23})x_2              + d_3x_4], \\
\dot{x_4} & = x_4[ -1+\mathcal{X} - c_{14}x_1 + (d_1(-1+c_{14})-t_{24})x_2 + e_{34}x_3 - d_3(1+e_{34})x_4],
\end{aligned}
\end{equation}
where $\mathcal{X}=x_1+x_2+x_3+x_4$. 
The coefficients are written in terms of the eigenvalues, classified as contracting, expanding and transverse. 
We denote by $-c_{ij}$ the contracting eigenvalue at $\xi_i$ in the direction of the $j^{th}$ basis vector, and analogously by $e_{ij}$ and by $-t_{ij}$ the expanding and transverse eigenvalues, respectively. 
The radial eigenvalues at $\xi_1$ and $\xi_3$ are both~$-1$.
At $\xi_2$ and $\xi_4$, the radial eigenvalues are eigenvalues of the two matrices
$$
\left(
\begin{array}{cc}
-d_1      & d_1(d_1-1) \\
1+e_{12}  & 1 - d_1(1+e_{12})
\end{array}
\right)
\quad\text{and}\quad
\left(
\begin{array}{cc}
-d_3      & d_3(d_3-1) \\
1+e_{34}  & 1 - d_3(1+e_{34})
\end{array}
\right).
$$
Stability in the radial direction at~$\xi_2$ and~$\xi_4$ can be achieved by requiring $d_1>1/(2+e_{12})>\frac12$ and $d_3>1/(2+e_{34})>\frac12$.
We note that complex radial eigenvalues are possible.


\begin{figure}
\hbox to \textwidth{\hfil\fbox{%
\includegraphics{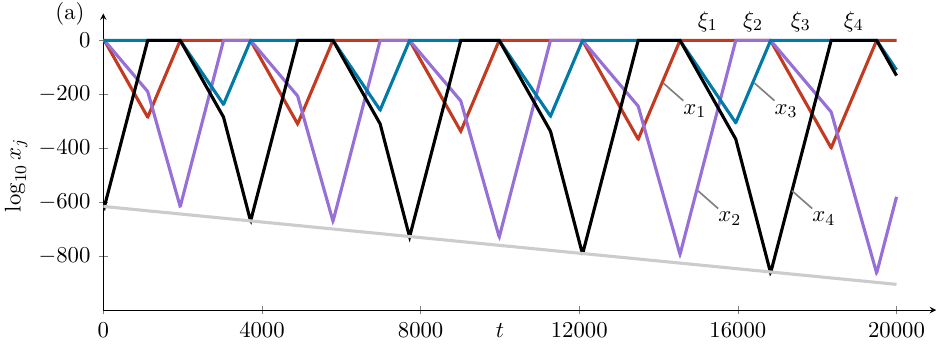}}\hfil}%

\hbox to \textwidth{\hfil\fbox{%
\includegraphics{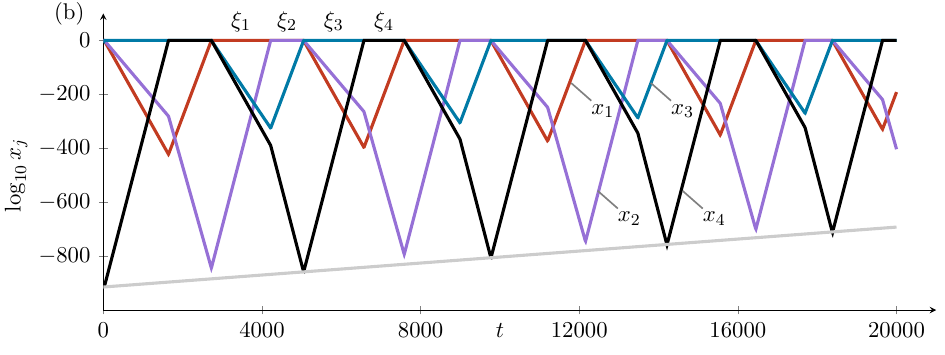}}\hfil}%

\caption{Illustration of the dynamics of Case~1~\eqref{eq:ODE-example-1}.
The parameters are $d_1=1.1$, $d_3=1.1$, $e_{12}=1.3$, $c_{13}=0.5$, $c_{14}=0.6$, $t_{24}=1.3$, $e_{34}=1.3$, $c_{31}=0.6$, $c_{32}=0.4$, $t_{42}=1.2$,
and (a) $e_{23}=0.8$, $e_{41}=0.8$ ($\delta_1=1.08654$, stable heteroclinic cycle);
    (b) $e_{23}=0.9$, $e_{41}=0.9$ ($\delta_1=0.93886$, unstable heteroclinic cycle).
The colors represent $x_1$~(red), $x_2$~(purple), $x_3$~(blue) and $x_4$~(black), plotted in logarithmic (base~10) coordinates. 
The initial conditions are $(x_1,x_2,x_3)=(1,d_1,10^{-10})$ and (a)~$x_4=10^{-600}$ and (b)~$x_4=10^{-900}$.
The grey line indicates a factor of~$\delta_1$ growth of the minima of~$\log x_4$ as well as a factor of~$\delta_1$ growth of the time interval between these minima, as indicated by~\eqref{eq:case1transition}.
Here, the grey line matches the successive minima of~$\log x_4$ (black), with $\log x_4$ growing in magnitude in the stable case and decreasing in magnitude in the unstable case.
In the unstable case, the trajectory eventually leaves the neighborhood of the heteroclinic cycle.
\label{fig:example1}}
\end{figure}

In Figure~\ref{fig:example1} we give examples of parameter values where the heteroclinic cycle is (a)~stable and (b)~unstable.
This example is based on the convection problem examined by~\cite{Matthews1996} in~$\R^7$, and it is capable of the same global bifurcations and chaotic dynamics reported in~\cite{Matthews1996}.

\subsection{Case 2}

The second example is a cycle between the four equilibria $\xi_1=(1,0,0,0)$, $\xi_2=(d_1,1,0,0)$, $\xi_3=(0,0,1,0)$ and $\xi_4=(0,0,0,1)$, where $d_1>0$.
As before, we require radial stability of the four equilibria, and set the eigenvalues of the origin in the $x_2$, $x_3$ and $x_4$ directions to be $-1$, $+1$, $+1$, equal in magnitude to simplify the presentation.
Expressing the coefficients in the ODE in terms of the eigenvalues, we have
\begin{equation}\label{eq:ODE-example-2}
\begin{aligned}
\dot{x_1} & = x_1[\phantom{-}1-\mathcal{X}             + d_1x_2                     - c_{31}x_3 + e_{41}x_4], \\
\dot{x_2} & = x_2[ -1+\mathcal{X} + e_{12}x_1 - d_1(1+e_{12})x_2           - c_{32}x_3 -t_{42}x_4], \\
\dot{x_3} & = x_3[\phantom{-}1-\mathcal{X} - t_{13}x_1 + (d_1(1+t_{13})+e_{23})x_2              - c_{43}x_4], \\
\dot{x_4} & = x_4[\phantom{-}1-\mathcal{X} - c_{14}x_1 + (d_1(1+c_{14})-t_{24})x_2 + e_{34}x_3],
\end{aligned}
\end{equation}
where $\mathcal{X}=x_1+x_2+x_3+x_4$. 
The radial eigenvalues at $\xi_1$, $\xi_3$ and $\xi_4$ are all~$-1$.
At $\xi_2$, the radial eigenvalues are eigenvalues of the matrix
$$
\left(
\begin{array}{cc}
-d_1      & d_1(d_1-1) \\
1+e_{12}  & 1 - d_1(1+e_{12})
\end{array}
\right).
$$
Stability in the radial direction at~$\xi_2$ can be achieved by requiring $d_1>1/(2+e_{12})>\frac12$.
We note that complex radial eigenvalues are possible.


\begin{figure}
\hbox to \textwidth{\hfil\fbox{%
\includegraphics{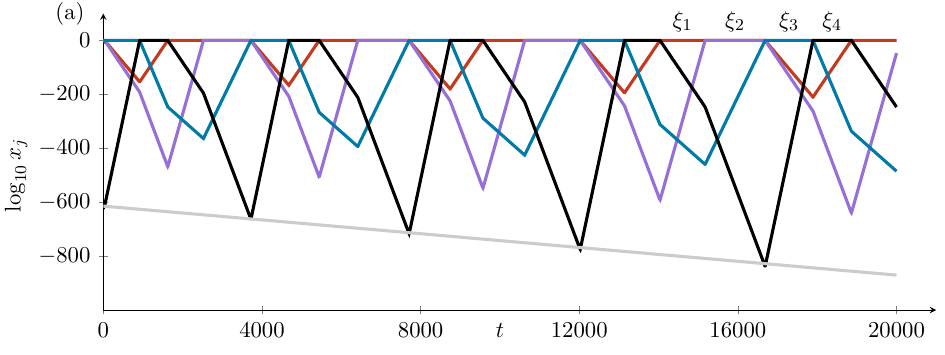}}\hfil}%

\hbox to \textwidth{\hfil\fbox{%
\includegraphics{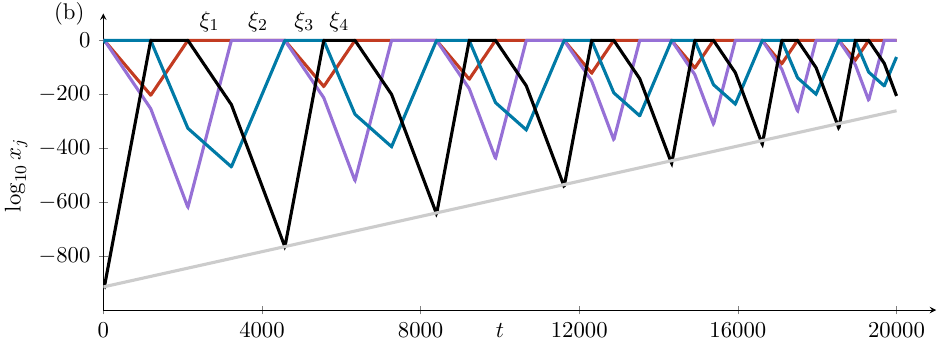}}\hfil}%

\caption{Illustration of the dynamics of Case~2~\eqref{eq:ODE-example-2}.
The parameters are $d_1=1.1$, $t_{13}=0.3$, $c_{14}=0.5$, $t_{24}=0.9$, $c_{31}=0.4$, $c_{32}=0.5$, $e_{41}=0.5$, $t_{42}=0.9$, $c_{43}=0.8$,
and (a)~$e_{12}=1.2$, $e_{23}=0.7$, $e_{34}=1.6$ ($\delta_2=1.07708$);
    (b) $e_{12}=1.3$, $e_{23}=0.8$, $e_{34}=1.8$ ($\delta_2=0.83665$).
The colors and initial conditions are as in Figure~\ref{fig:example1}.
The grey line indicates a factor of~$\delta_2$~\eqref{eq:case2transition} growth, which lines up well with the minima of~$\log x_4$ (black).
\label{fig:example2}}
\end{figure}

In Figure~\ref{fig:example2} we give examples of parameter values where the heteroclinic cycle is (a)~stable and (b)~unstable.

\subsection{Case 3}

The third example is a cycle between the four equilibria $\xi_1=(1,0,0,0)$, $\xi_2=(d_1,1,0,0)$, $\xi_3=(0,d_2,1,0)$ and $\xi_4=(0,0,0,1)$, where $d_1>0$ and $d_2>0$.
As before, we require radial stability of the four equilibria, and set 
the eigenvalues of the origin in the $x_2$, $x_3$ and $x_4$ directions to be $-1$, $+1$, $+1$, equal in magnitude to simplify the presentation.
Expressing the coefficients in the ODE in terms of the eigenvalues, we have
\begin{equation}\label{eq:ODE-example-3}
\begin{aligned}
\dot{x_1} & = x_1[\phantom{-}1-\mathcal{X}             + d_1x_2                     +(d_2(1-d_1)-c_{31})x_3 + e_{41}x_4], \\
\dot{x_2} & = x_2[ -1+\mathcal{X} + e_{12}x_1 - d_1(1+e_{12})x_2           + d_2(d_1(1+e_{12})-1)x_3 -c_{42}x_4], \\
\dot{x_3} & = x_3[\phantom{-}1-\mathcal{X} - t_{13}x_1 + (d_1(1+t_{13})+e_{23})x_2 + d_2(1-d_1(1+t_{13})-e_{23})x_3              - c_{43}x_4], \\
\dot{x_4} & = x_4[\phantom{-}1-\mathcal{X} - c_{14}x_1 + (d_1(1+c_{14})-t_{24})x_2 + (d_2(1+t_{24})-d_1d_2(1+c_{14})+e_{34})x_3],
\end{aligned}
\end{equation}
where $\mathcal{X}=x_1+x_2+x_3+x_4$. 

The radial eigenvalues at $\xi_1$ and $\xi_4$ are both~$-1$.
At $\xi_2$ and~$\xi_3$, the radial eigenvalues are eigenvalues of the two matrices
$$
\left(\!
\begin{array}{cc}
-d_1      & d_1(d_1-1) \\
1+e_{12}  & 1 - d_1(1+e_{12})
\end{array}
\!\right)
\,\text{and}\,
\left(\!
\begin{array}{cc}
d_2(1-d_1(1+e_{12}))     & d_2(1-d_2+d_1d_2(1+e_{12})) \\
d_1(1+t_{13})-1+e_{23}  & d_2(1-e_{23}) - d_1d_2(1+t_{13}) - 1
\end{array}
\!\right).
$$
Stability in the radial direction at~$\xi_2$ can be achieved by requiring $d_1>1/(2+e_{12})>\frac12$. 
Radial stability at~$\xi_3$ is more complicated but can be readily checked in individual examples.
We note that complex radial eigenvalues are possible.


\begin{figure}
\hbox to \textwidth{\hfil\fbox{%
\includegraphics{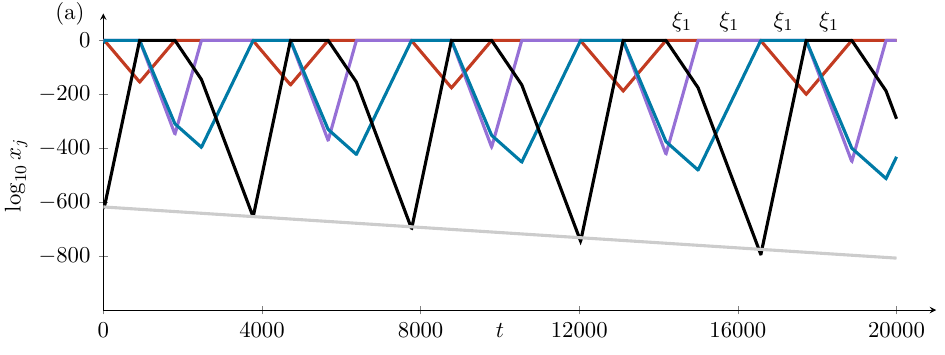}}\hfil}%

\hbox to \textwidth{\hfil\fbox{%
\includegraphics{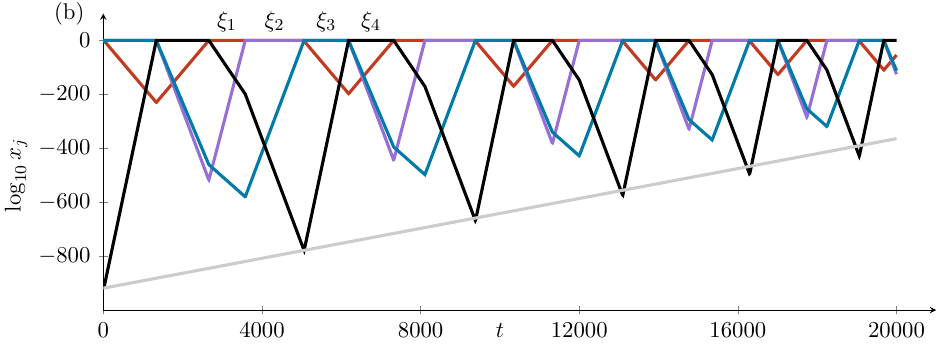}}\hfil}%

\caption{Illustration of the dynamics of Case~3~\eqref{eq:ODE-example-3}.
The parameters are $d_1=1.1$, $d_2=1$, $t_{13}=0.3$, $c_{14}=0.5$, $t_{24}=0.9$, $e_{34}=1.6$, $c_{31}=0.4$, $e_{41}=0.4$, $c_{42}=0.9$, $c_{43}=0.8$,
and (a)~$e_{12}=1.2$, $e_{23}=0.7$, ($\delta_3=1.05804$);
    (b)~$e_{12}=1.3$, $e_{23}=0.9$, ($\delta_3=0.84615$).
The colors and initial conditions are as in Figure~\ref{fig:example1}.
The grey line indicates a factor of~$\delta_3$~\eqref{eq:case3transition} growth, which lines up reasonably well with the minima of~$\log x_4$ (black).
\label{fig:example3}}
\end{figure}

In Figure~\ref{fig:example3} we give examples of parameter values where the heteroclinic cycle is (a)~stable and (b)~unstable.

\subsection{Case 4}

The fourth example is a cycle between the four equilibria $\xi_1=(1,0,0,0)$, $\xi_2=(d_1,1,0,0)$, $\xi_3=(0,d_2,1,0)$ and $\xi_4=(0,0,d_3,1)$, where $d_1>0$, $d_2>0$ and $d_3>0$.
As before, we require radial stability of the four equilibria, and it turns out that setting 
the eigenvalues of the origin in the $x_2$, $x_3$ and $x_4$ directions to be $-1$, $+1$, $+1$ is helpful for this.
Expressing the coefficients in the ODE in terms of the eigenvalues, we have
\begin{equation}\label{eq:ODE-example-4}
\begin{aligned}
\dot{x_1} & = x_1[\phantom{-}1-\mathcal{X}             + d_1x_2                     +(d_2(1-d_1)-c_{31})x_3 + (d_3(1+c_{31})+d_2d_3(d_1-1)+e_{41})x_4], \\
\dot{x_2} & = x_2[ -1+\mathcal{X} + e_{12}x_1 - d_1(1+e_{12})x_2           + d_2(d_1(1+e_{12})-1)x_3 \\
 & \qquad\qquad\qquad\qquad\qquad\qquad{}-(d_1d_2d_3(1+e_{12})+d_3(1-d_2)+c_{42})x_4], \\
\dot{x_3} & = x_3[\phantom{-}1-\mathcal{X} - c_{13}x_1 + (d_1(1+c_{13})+e_{23})x_2 + d_2(1-d_1(1+c_{13})-e_{23})x_3 \\
 & \qquad\qquad\qquad\qquad\qquad\qquad{}+ d_3(d_1d_2(1+c_{13}) + d_2(e_{23}-1)+1)x_4], \\
\dot{x_4} & = x_4[\phantom{-}1-\mathcal{X} - c_{14}x_1 + (d_1(1+c_{14})-t_{24})x_2 + (d_2(1+t_{24})-d_1d_2(1+c_{14})+e_{34})x_3 \\
 & \qquad\qquad\qquad\qquad\qquad\qquad{}+ d_3(d_1d_2(1+c_{14})-d_2(1+t_{24})-e_{34}+1)x_4],
\end{aligned}
\end{equation}
where $\mathcal{X}=x_1+x_2+x_3+x_4$. 

The radial eigenvalue at $\xi_1$ is~$-1$.
At $\xi_2$, $\xi_3$ and~$\xi_4$, the radial eigenvalues are eigenvalues of three $2\times2$ matrices. The first two of these are the same as in Example~3 (apart from relabelling $t_{13}$ as $c_{13}$); the third (for~$\xi_4$) is
$$
\left(\!
\begin{array}{cc}
-d_3(d_1d_2(1+c_{13})+d_2(e_{23}-1)+1)      & 
    d_3(d_1d_2d_3(1+c_{13})+d_2d_3(e_{23}-1)+d_3-1) \\
d_2(1+t_{24})-d_1d_2(1+c_{14})+e_{34}-1  & 
d_1d_2d_3(1+c_{14})-d_2d_3(1+t_{24})+d_3(1-e_{34})-1
\end{array}
\!\right).
$$
Radial stability can readily be checked in individual examples.


\begin{figure}
\hbox to \textwidth{\hfil\fbox{%
\includegraphics{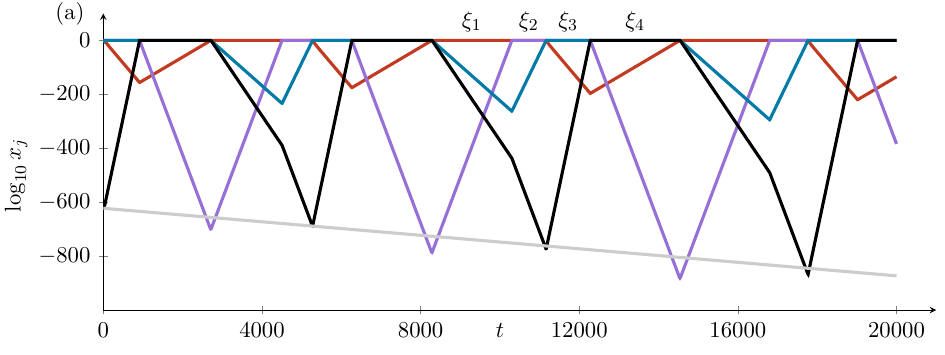}}\hfil}%

\hbox to \textwidth{\hfil\fbox{%
\includegraphics{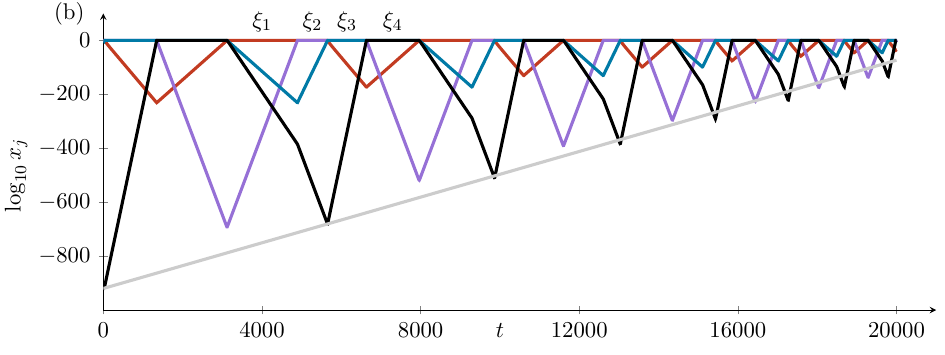}}\hfil}%

\caption{Illustration of the dynamics of Case~4~\eqref{eq:ODE-example-4}.
The parameters are $d_1=1.3$, $d_2=1.1$, $d_3=1.1$, $e_{12}=0.9$, $c_{13}=0.3$, $c_{14}=0.5$, $e_{23}=0.7$, $t_{24}=0.9$, $e_{34}=1.6$, $c_{31}=0.4$, $c_{42}=0.9$,
and (a)~$e_{41}=0.2$ ($\delta_4=1.10714$);
    (b)~$e_{41}=0.3$ ($\delta_4=0.73810$).
The colors and initial conditions are as in Figure~\ref{fig:example1}.
The grey line indicates a factor of~$\delta_4$~\eqref{eq:case4transition} growth, which lines up approximately with the minima of~$\log x_4$ (black), deviating in particular in the unstable case as the trajectory leaves the heteroclinic cycle.
\label{fig:example4}}
\end{figure}

In Figure~\ref{fig:example4} we give examples of parameter values where the heteroclinic cycle is (a)~stable and (b)~unstable. In this case, the grey line, indicting a factor of $\delta_4$ growth in successive time intervals and minima of~$\log x_4$, is noticeably different from the actual locations of the minima. Decreasing the initial conditions from $x_4=10^{-600}$ and $x_4=10^{-900}$ to (for example) $x_4=10^{-6000}$ reduces this discrepancy, while for an initial condition of $x_4=10^{-60}$, the discrepancy is even more pronounced. This observation applies to the other cases as well: the discrepancy arises from the $\Order(\log A)$ and $\Order(\log h)$ terms that have been dropped in deriving \eqref{eq:case1transition}--\eqref{eq:case4transition}.

\section{Discussion}\label{sec:discussion}

Our results provide a starting point for a general approach to the study of the stability of a broader class of robust heteroclinic cycles.
Up until now the systematic approaches to stability required the existence of contracting eigenvalues at every equilibrium.
We have shown how to treat the absence of contracting eigenvalues, and although our specific examples are in $\R^4$, the principles of the calculations are applicable to any dimension.
Each transition map from $H_{i}^{in,i-1}$ to $H_{i+1}^{in,i}$ depends on the locations of~$\xi_{i-1}$, $\xi_i$ and~$\xi_{i+1}$ in these examples.
In higher dimensions, there will be a greater variety of possible transition matrices, with many possible combinations of dimensions.

Another interesting feature of robust heteroclinic cycles in pluridimensions is that there are also equilibria with more than one contracting directions. 
We have shown that, on the incoming section, the values of the contracting coordinates do not contribute to the stability calculation, though the contracting eigenvalues do. 

We have presented all examples of cycles in pluridimensions in~$\R^4$ satisfying (A1) (one-dimensional unstable manifolds), (A2) (invariant coordinate axes and hyperplanes), (A3) (one equilibrium on each connected component) and (A4) (the origin is excluded).
These assumptions were only introduced for the purpose of constructing simple examples but are not required for 
a robust heteroclinic cycle to have $P$ subspaces that vary in dimension around the cycle.
Examples of robust cycles in pluridimensions not satisfying these assumptions can be treated in a similar manner, with different degrees of extra effort.
In turn:
\begin{itemize}
\item
Allowing higher-dimensional unstable manifolds, going beyond Definition~\ref{def:simple-HS} and relaxing Assumption~(A1), would bring in aspects of cycles with two (or more) dimensional connections, as in~\cite{Castro2022b}, or heteroclinic networks as in (for example)~\cite{Kirk1994}.
The dynamics near such a network can involve trajectories making choices as to which direction to take and the stability of trajectories following sequences of choices is already understood in terms of transition matrices \cite{Krupa2004,Podvigina2012,Podvigina2023,Postlethwaite2022}.
It would be very interesting to bring the systematic handling of $P$ subspaces that are of different dimension to the theory of heteroclinic networks.

\item
Our examples do not require any symmetries, but our assumption (A2) about invariant coordinate axes and hyperplanes could be replaced, using symmetries to guarantee the structurally stable connections needed for a heteroclinic cycle.
Symmetries, for example having reflection symmetry in every coordinate, can lead to all hyperplanes being invariant in the same way as~(A2).
However, symmetries can act more generally than this.
In the example of~\cite{Matthews1996} in~$\R^7$, many of the coordinate planes are not invariant, and the sequence of dimensions of the $L$~subspaces is $2\rightarrow4\rightarrow2\rightarrow4$, and that of the $P$~subspaces is $4\rightarrow5\rightarrow4\rightarrow5$.
So, although the unstable manifolds of the equilibria are still one-dimensional, we do not always have~(C2): there are equilibria with $\dim P_i>\dim L_i+1$.
This happens because some of the equilibria have negative as well as positive expanding eigenvalues: this is prevented by~(A2).
In the related model of~\cite{Rucklidge1995b} in $\R^9$, some of the variables can change sign as they approach the heteroclinic cycle.
This is also prevented by~(A2).
Even so, we expect that the approach to calculating stability that we have taken here will work, with appropriate modifications, in these two examples.

\item
The example of~\cite{Sikder1994} does not satisfy Assumption~(A3) but stability calculations for that, and similar, examples would carry through unchanged.
Allowing more than one equilibrium on an axis within a cycle could lead to further interesting generalisations~\cite{Ashwin2013,Castro2023}.

\item
The example of~\cite{Hawker2005a} in~$\R^3$ includes the origin (as well as two on-axis equilibria) and so does not satisfy~(A4).
The origin has two contracting directions, the first of the on-axis equilibria has a transverse but no contracting direction, and the second on-axis equilibrium has one contracting direction.
The stability calculations can be handled in a similar way, and we find that the stability of the cycle is determined by the product of three transition matrices, of the form of~\eqref{eq:mapPPX}, \eqref{eq:mapPAX} and~\eqref{eq:mapAPX}. 
When multiplied out, our method agrees with the results of~\cite{Hawker2005a}.
Just as in Section~\ref{sec:transitionmatrices} in the plane-to-axis-to-plane case, the values of the logarithms of the two contracting variables at the origin are both $\Order(\log h)$ and so they can both be neglected.
\end{itemize}

We remark that heteroclinic cycles in systems with symmetry are often associated with certain patterns in the lattice of isotropy subgroups~\cite{Melbourne1989b}, where equilibria in maximal fixed point subspaces are linked by connections in submaximal fixed point subspaces.
Up-and-down patterns in the lattice of isotropy subgroups indicate the possibility of robust heteroclinic cycles.
In the examples in~\cite{Matthews1996,Rucklidge1995b}, the pluridimensional nature of the examples is related to the fact that the connections \emph{skip} a level in the lattice of isotropy subgroups.
This observation suggest that heteroclinic cycles in pluridimensions might be sought in symmetric systems having lattices of isotropy subgroups with sufficiently many levels: this will be a subject of future work.

We end by observing that robust cycles in pluridimensions 
form an important class of non-simple heteroclinic cycles, and the work we have presented here is a starting point to a general theory of their stability.
This type of cycle will arise, for example, in modelling the dynamics of evolving populations when there are transitions between equilibria corresponding to mixed populations with different numbers of species, as in the example of~\cite{Matias2018}.
We also expect that these ideas we have presented will be useful for analysing other more general problems, such as the stability of depth~2 heteroclinic cycles~\cite{Chawanya1997}.

\section*{Acknowledgments}
We would like to acknowledge conversations with Alexander Lohse and 
Claire Postlethwaite, as well as useful remarks from Peter Ashwin, Christian Bick, Martin Golubitsky and Josef Hofbauer.

The first author was partially supported by CMUP, member of LASI, which is financed by national funds through FCT -- Funda\c{c}\~ao para a Ci\^encia e a Tecnologia, I.P., under the projects UIDB/00144/2020 and UIDP/00144/2020.
This work started when the first author visited the University of Leeds, whose hospitality and financial support through the EPSRC grant EP/V014439/1 is gratefully acknowledged.
We are also grateful for the support of the \hbox{CMUP} and the hospitality of the University of Porto during a visit of the second author to Porto.
The data associated with this paper are openly available from the University of Leeds Data Repository (\url{http://doi.org/10.5518/1494})~\cite{Castro2024b}, as is the program that generated the data.
For the purpose of open access, the authors have applied a Creative Commons Attribution (CC~BY) license to any Author Accepted Manuscript version arising from this submission.


\end{document}